\documentclass[12pt]{article}
\usepackage{amsfonts}
\usepackage{amsmath}
\usepackage{array}
\usepackage{graphicx}
\usepackage{xcolor}
\usepackage{float}
\usepackage{appendix}
\usepackage{subfig}
\usepackage{longtable}
\usepackage{rotating}
\usepackage{framed}
\usepackage{algorithm}
\usepackage{algorithmic}
\usepackage{pdflscape}
\usepackage{afterpage}
\usepackage{booktabs}
\usepackage{caption}
\usepackage{hyperref}
\usepackage{appendix}
\usepackage{biblatex}
\usepackage{comment}

\hypersetup{
	colorlinks=true,
	linkcolor=blue,
	filecolor=magenta,      
	urlcolor=cyan,
}
\newtheorem{theorem}{Theorem}

\newtheorem{corollary}[theorem]{Corollary}

\newtheorem{definition}{Definition}
\newtheorem{example}[theorem]{Example}

\newtheorem{proposition}[theorem]{Proposition}
\newtheorem{remark}{Remark}

\newenvironment{proof}[1][Proof]{\noindent\textbf{#1.} }{\ \rule{0.5em}{0.5em}}
\input{tcilatex}
\begin{document}

\title{\textbf{Robust and efficient Breusch-Pagan test-statistic: an
application of the }$\beta $\textbf{-score Lagrange multipliers test for
non-identically distributed individuals }}
\date{\today}
\author{Nirian Mart\'{\i}n, Complutense University of Madrid (\texttt{\href{mailto:nimartin@ucm.es}%
		{nimartin@ucm.es}})}
\maketitle
\begin{abstract}
In Econometrics, the Breusch-Pagan test-statistic has become an iconic
application of the Lagrange multipliers (LM) test. We
shall introduce $\beta $-score LM tests for heteroscedasticity in linear
regression models, which trades-off the degree of robustness and efficiency is
through a tuning parameter $\beta \geq 0$, being $\beta =0$ the
classical Breusch-Pagan test-statistic, the most efficient one under absence
of outliers. A very elegant expression is obtained, with an appealing least squares interpretation. The construction of the test-statistic is performed extending the
methodology of Basu et al. (2022) from identically distributed\ to
non-identically distributed individuals, for composite null hypotheses.
Detailed theoretical justifications about robustness and efficiency
properties are given, all of them under normality. A modified version is derived, the Koenker's $\beta $-score test-statistic.
The way of dealing, in practice, with new heteroscedasticity tests for linear
regression models is shown through a classical example.
\end{abstract}

\noindent \textbf{JEL CLASSIFICATION}: C01, C10;

\medskip

\section{Introduction}

Majority of undergraduate textbooks of Econometrics, include the Breusch-Pagan
Lagrange Multipliers (LM) test as one of the most important tests for
heteroscedasticity in linear regression models. Authored by Breusch and Pagan
(1979), it was also obtained independently by Godfrey (1978a) under a stronger
assumption of multiplicative heteroskedasticity. Thereafter, Koenker (1981)
and Koenker and Basset (1982) developed a modified version of the test that
relaxes the assumption of normality in the error distribution, allowing for a
more general distribution. With a particular structure of the design matrix
for the type of heteroskedasticity, adding the cross-product of the
regressors, a test for heteroscedasticity was obtained in White (1980), but
Waldman (1983) clarified that White's version of the test is a particular case
of the Breusch-Pagan or Koenker's tests, depending on the normal ot non-normal
distributional assumption for errors. The mentioned tests suffer from lack of
robustness under outlying data. This issue is reported for example in Carrol
and Ruppert (1988, page 98), Green (page 314), Kmenta (1986, page 295), Lyon
and Tsay (1996, page 339) and Kalina (2011). Recently, Berenguer-Rico and
Wilms (2021) have suggested outlier removal for testing heterogeneity using
the White test-statistic, while Alih and Ong (2015) have proposed a robust
version of the Goldfeld-Quandt test-statistic replacing its non-robust component.

The robustness of many statistical procedures is often achieved at the price
of the loss of efficiency. The estimators and test-statistics constructed with
simultaneous asymptotic efficiency and robustness properties are very
attractive. This idea was firstly introduced, in the estimation setting, in
Beran (1977) through the Hellinger distance. The original notion of distance
(or divergence) between two probability distributions, constructed from a
sample, was introduced separately by Kolmogorov (1933) and Mahalanobis (1936),
in different ways. In parametric statistical inference, from the beginning and
in successive developments, the proposed new statistical distances were
focussed on the closeness of a empirical distribution $g$ and a model based
distribution $f_{\boldsymbol{\theta}}$, where $\boldsymbol{\theta}$\ is the
$p$-dimensional parameter vector of interest. When the true distribution
belongs to the model based one, it is well-known that the Kullback-Leibler
(1951) divergence covers all the classical statistical theory related to the
maximum likelihood estimators (MLEs) and their associated test-statistics.
There is an extent family of $\phi$-divergence measures, introduced by
Csisz\'{a}r (1967), from which BAN (Best Asymptotically Normal) estimators are
obtained, i.e. estimators as efficient as the MLEs, asymptotically. The
aforementioned distances include as a particular class of distances, the so
called power divergences\ of Cressie and Read (1984) and the last ones at the
same time contain the Hellinger distance as a particular member. For more
details see Pardo (2006) and Basu et al. (2011). With the progress of the
literature, several new robust and efficient\ minimum distance procedures have
been proposed and among them, the density power divergence (DPD) measures of
Basu et al. (1998)\ have become very popular. Their novelty with respect to
the Hellinger distance was based on the bounded influence function of the
estimators as well as on the extension of their validity to continuous
populations apart from the discrete ones, mostly treated until this moment.
Thereafter, the DPDs have received a growing attention in statistical
inference being applied for estimation and testing in different parametric
models. For linear regression, from a conditional point of view or fixed
design matrix, Ghosh and Basu (2013) established for the first time formally
how to deal with DPDs in estimation when the observations in the sample are
independent and not-identically distributed. For robust testing through DPDs
in linear regression models, in Section 6.3. of Ghosh et al. (2016) $\beta
$-Wald test statistics were proposed, while in Qin and Priebe (2017) the
$\beta$-likelihood ratio test statistics (L$q$-likelihood-ratio-type tests)
were introduced. There is a sounded disadvantage of the $\beta$-Likelihood
Ratio test statistics with respect to the $\beta$-Wald test statistics, the
asymptotic distribution depends on weights to be calculated as eigenvalues of
matrices dependent of unknown parameters and except for scalar parameters this
issue affects the accuracy of the calculation of the test-statistics. Based on
Basu et al. (2022), it is expected that the $\beta$-score LMs to be developed
in the current paper, will have the same asymptotic distribution as the
$\beta$-Wald test statistics. In addition, as noted by Boos (1992), the
strength of the generalized score tests, where the $\beta$-score LMs are
included, is the invariance property. This is in fact, what happens with the
Breusch-Pagan LM test, since it is valid for testing the presence of a general
type of heteroscedasticity. A unique expression of the test-statistic covers
different possible expressions of scedasticity function, which is not possible
with other types of test-statistics.

In econometrics, there have been substantial contributions to assess
robustness while retaining little efficiency, based on minimum DPD estimators
(for example, see Lee and Song (2009), Kim and Lee (2013, 2017, 2018)). To the best of our knowledge, there
are no publications about testing using such estimators, except for goodness-of-fit setting (for example, Kim
(2018)).
Based on minimum DPD estimators, consistent estimators are considered, as well
as unbiased estimating equations. The proposal of this paper is proven to be
theoretically robust two-fold, in constructing the estimators and also the
test-statistic. The main feature of the methodology is on one hand in
achieving robustness through a bounded influence function for estimation and
testing and with an increasing gross error sensitivity as the tuning parameter
$\beta$ increases, while also ensuring an acceptable level of efficiency
through the same tuning parameter $\beta$ approaching $0$. On the other hand,
the Pitman's asymptotic relative efficiency (ARE) is obtained, which measures
the price, in terms of the relative sample size, to get a particular value of
the asymptotic power when dealing with pure data.

The paper is organized as follows. Based on a new framework of sampling, with
independent but non-identically observations, for composite hypothesis testing
is introduced and motivated in Section \ref{Sec1}, by relating it to the
heteroscedastic linear regression model. The main theoretical results are
presented in Section \ref{Sec2}, derived the Breusch-Pagan $\beta$-score LM
tests in Section \ref{Sec3} and influence function analysis in Section
\ref{Sec4}. A modified version is derived in Section \ref{Sec4b}, the Koenker's $\beta $-score test-statistic. 
In Section \ref{Sec5}, a practical demonstration of how to handle new tests for heterogeneity in linear regression models is presented using a well-known example. Some concluding remarks are given in Section
\ref{Sec6}.
Most of the proofs are collected in the Appendix, at the end of the paper.

\section{Basic model and testing specifications\label{Sec1}}

Based on the conditional version of the linear regression, the predictors are considered fixed, so the sample of responses $Y_{i}$, $i=1,\ldots
,n$, requires adapting the previous theory presented in Basu et al. (2002) for non-identically distributed individuals. In this setting, being
$g_{i}(\cdot)$ the true probability density function (p.d.f.) for
$i=1,\ldots,n$, and $f_{i,\boldsymbol{\theta}}(\cdot)$ the p.d.f. under the
model with $\boldsymbol{\theta}\in\Theta\subset%
\mathbb{R}
^{p}$, the density power divergence of the whole sample is, according to Ghosh
and Basu (2013), given by $H_{n,\beta}(\boldsymbol{\theta})=\sum_{i=1}%
^{n}d_{\beta}(g_{i},f_{i,\boldsymbol{\theta}})$, where%
\[
d_{\beta}(g_{i},f_{i,\boldsymbol{\theta}})=\int_{\mathcal{Y}}\left(
f_{i,\boldsymbol{\theta}}^{\beta+1}(y)-\tfrac{\beta+1}{\beta}%
f_{i,\boldsymbol{\theta}}^{\beta}(y)g_{i}+\tfrac{1}{\beta}g_{i}^{\beta
	+1}\right)  dy,\quad\beta>0.
\]
where $\mathcal{Y}$ is the common support for the sample of all the responses.
In practice, the unknown p.d.f. $g_{i}(\cdot)$\ is approximated through the
empirical one, $\widehat{g}_{i}(\cdot)$.

$H_{0}:\boldsymbol{\theta }\in \Theta _{0}$\ vs.$\ H_{1}:%
\boldsymbol{\theta }\in \Theta -\Theta _{0}$, with 
\begin{equation}
\Theta _{0}=\{\boldsymbol{\theta }\in \Theta :\boldsymbol{m}(\boldsymbol{%
\theta })=\boldsymbol{0}_{r}\},  \label{eq4}
\end{equation}%
where $\boldsymbol{m}:\Theta \longrightarrow 
\mathbb{R}
^{r}$, with $r<p$ is a composite null hypothesis test. It is assumed the
function which defines $\Theta _{0}$ to fulfill some regularity conditions,%
\begin{equation*}
\boldsymbol{M}(\boldsymbol{\theta })=\frac{\partial }{\partial \boldsymbol{%
\theta }}\boldsymbol{m}^{T}(\boldsymbol{\theta }),
\end{equation*}%
exists and is continuous in $\boldsymbol{\theta }$ and $\mathrm{rank}(%
\boldsymbol{M}(\boldsymbol{\theta }))=r$.

\begin{definition}
	The minimum DPD estimator of $\boldsymbol{\theta}$, restricted to null
	hypothesis established by (\ref{eq4}), is obtained as
	\[
	\widetilde{\boldsymbol{\theta}}_{n,\beta}=\arg\min_{\boldsymbol{\theta}%
		\in\Theta_{0}}\widehat{H}_{n,\beta}(\boldsymbol{\theta})=\arg\max
	_{\boldsymbol{\theta}\in\Theta_{0}}\left(  \tfrac{\beta+1}{\beta}\tfrac{1}%
	{n}\sum_{i=1}^{n}f_{i,\boldsymbol{\theta}}^{\beta}(Y_{i})-\sum_{i=1}%
	^{n}\mathrm{E}\left[  f_{i,\boldsymbol{\theta}}^{\beta}(Y_{i})\right]
	\right)  ,
	\]
	where
	\begin{align*}
		\widehat{H}_{n,\beta}(\boldsymbol{\theta}) &  =\sum_{i=1}^{n}d_{\beta
		}(\widehat{g}_{i},f_{i,\boldsymbol{\theta}}),\\
		d_{\beta}(\widehat{g}_{i},f_{i,\boldsymbol{\theta}}) &  =\mathrm{E}\left[
		f_{i,\boldsymbol{\theta}}^{\beta}(Y_{i})\right]  -\tfrac{\beta+1}{\beta}%
		\tfrac{1}{n}f_{i,\boldsymbol{\theta}}^{\beta}(Y_{i})+\tfrac{1}{\beta}\tfrac
		{1}{n}.
	\end{align*}
	
\end{definition}

The restricted maximum likelihood estimator (MLE) of $\boldsymbol{\theta }$
is a member of the minimum DPD estimator, since%
\begin{equation*}
\widetilde{\boldsymbol{\theta }}_{\beta =0}=\arg \min_{\boldsymbol{\theta }%
\in \Theta _{0}}\lim_{\beta \rightarrow 0^{+}}\widehat{H}_{n,\beta }(%
\boldsymbol{\theta })=\arg \max_{\boldsymbol{\theta }\in \Theta
_{0}}\sum_{i=1}^{n}\log f_{i,\boldsymbol{\theta }}(Y_{i})=\widetilde{%
\boldsymbol{\theta }}.
\end{equation*}

\begin{proposition}
	The minimum DPD estimator of $\boldsymbol{\theta}$, restricted to null
	hypothesis established by (\ref{eq4}), is obtained as solution in
	$(\boldsymbol{\theta}^{T},\boldsymbol{\lambda}_{n,\beta}^{T})^{T}$ of the
	following system of $p+r$ equations and unknown parameters
	\begin{align}
		\boldsymbol{U}_{n,\beta}\left(  \boldsymbol{\theta}\right)  -\boldsymbol{M}%
		\left(  \boldsymbol{\theta}\right)  \boldsymbol{\lambda}_{n,\beta} &
		=\boldsymbol{0}_{p},\label{eqA}\\
		\boldsymbol{m}\left(  \boldsymbol{\theta}\right)   &  =\boldsymbol{0}%
		_{r},\label{eqB}%
	\end{align}
	where $\boldsymbol{\lambda}_{n,\beta}$ is an $r$-vector of Lagrange
	multipliers, whose minimum DPD estimator is denoted by $\boldsymbol{\lambda
	}_{n,\beta}(\widetilde{\boldsymbol{\theta}}_{\beta})$ and
	\[
	\boldsymbol{U}_{n,\beta}\left(  \boldsymbol{\theta}\right)  =\frac{1}{n}%
	\sum_{i=1}^{n}\boldsymbol{u}_{i,\beta}(Y_{i};\boldsymbol{\theta}),
	\]
	is the estimating function of the unrestricted minimum DPD estimator of
	$\boldsymbol{\theta}$, where
	\begin{align}
		\boldsymbol{u}_{i,\beta}(y;\boldsymbol{\theta}) &  =f_{i,\boldsymbol{\theta}%
		}^{\beta}(y)\boldsymbol{s}_{i,\boldsymbol{\theta}}(y)-\boldsymbol{\xi
		}_{i,\beta}\left(  \boldsymbol{\theta}\right)  ,\nonumber\\
		\boldsymbol{s}_{i,\boldsymbol{\theta}}(y) &  =\tfrac{\partial}{\partial
			\boldsymbol{\theta}}\log f_{i,\boldsymbol{\theta}}(y),\nonumber\\
		\boldsymbol{\xi}_{i,\beta}(\boldsymbol{\theta}) &  =\mathrm{E}\left[
		\boldsymbol{u}_{i,\beta}\left(  Y_{i},\boldsymbol{\theta}\right)  \right]
		=\int_{\mathcal{-\infty}}^{+\infty}f_{i,\boldsymbol{\theta}}^{\beta
			+1}(y)\boldsymbol{s}_{i,\boldsymbol{\theta}}(y)dy=\mathrm{E}%
		[f_{i,\boldsymbol{\theta}}^{\beta}(Y_{i})\boldsymbol{s}_{i,\boldsymbol{\theta
		}}(Y_{i})].\label{psi}%
	\end{align}
	
\end{proposition}

The so-called scedastic function, $h(\cdot )$, determines functionally the form of heteroscedasticity but it is not prefixed. It is assumed to be continuous,
to possess at least first and second derivarites and to verify 
\begin{equation*}
	h(\eta )>0\;\forall \eta \in 
	\mathbb{R}
	,\quad h(0)=1,\quad \left. \tfrac{d}{d\eta }h(\eta )\right\vert _{\eta
		=0}\neq 0.
\end{equation*}%
\begin{definition}\label{def0}
	The conditional heteroscedastic linear regression model is given by
	\begin{equation}
		Y_{i}=\boldsymbol{x}_{i}^{T}\boldsymbol{\beta}+\epsilon_{i},\text{\quad
		}i=1,\ldots,n,\label{eq1}%
	\end{equation}
	where $\epsilon_{i}\overset{ind}{\sim}\mathcal{N}(0,\sigma_{i}^{2})$ and
	\begin{align*}
		\boldsymbol{\beta} &  =(\beta_{0},\beta_{1},\ldots,\beta_{p})^{T}\in%
		\mathbb{R}
		^{p+1},\\
		\boldsymbol{x}_{i}^{T} &  =(1,x_{i1},\ldots,x_{ip})\in%
		\mathbb{R}
		^{p+1},\quad p\in%
		\mathbb{N}
		,\\
		\boldsymbol{\alpha} &  =(\alpha_{1},\ldots,\alpha_{r})^{T}\in%
		\mathbb{R}
		^{r},\\
		\boldsymbol{z}_{i}^{T} &  =(z_{i1},\ldots,z_{ir})\in%
		\mathbb{R}
		^{r},\quad r\in%
		\mathbb{N}
		,\\
		\sigma_{i}^{2} &  =\sigma_{i}^{2}(\boldsymbol{\alpha},\sigma^{2})=\sigma
		^{2}h(\boldsymbol{z}_{i}^{T}\boldsymbol{\alpha}).
	\end{align*}
	
\end{definition}
$h(\cdot )$ includes most
of the schemes considered in the literature, for example, the additive scedastic model $h(\boldsymbol{z}%
_{i}^{T}\boldsymbol{\alpha })=(1+\boldsymbol{z}_{i}^{T}\boldsymbol{\alpha }%
)$ or the multiplcative scedastic model $h(\boldsymbol{z}_{i}^{T}\boldsymbol{\alpha })=\exp \{\boldsymbol{z}%
_{i}^{T}\boldsymbol{\alpha }\}$. The explanatory variables, $\boldsymbol{x}%
_{i}$, $\boldsymbol{z}_{i}$, $i=1,\ldots ,n$, are assumed to be fixed, i.e.
non-random. In a full matrix notation, (\ref{eq1}) is given by%
\begin{equation}
\boldsymbol{Y}=\mathbb{X}_{n}\boldsymbol{\beta }+\boldsymbol{\epsilon },
\label{eq2}
\end{equation}%
where%
\begin{align*}
\boldsymbol{Y}& =(Y_{1},\ldots ,Y_{n})^{T}\sim \mathcal{N}(\mathbb{X}_{n}%
\boldsymbol{\beta },\mathrm{diag}\{\sigma _{i}^{2}(\boldsymbol{\alpha }%
,\sigma ^{2})\}_{i=1}^{n}), \\
\mathbb{X}_{n}& =(\boldsymbol{x}_{1},\ldots ,\boldsymbol{x}_{n})^{T}, \\
\boldsymbol{\epsilon }& =(\epsilon _{1},\ldots ,\epsilon _{n})^{T}\sim 
\mathcal{N}(\boldsymbol{0}_{n},\mathrm{diag}\{\sigma _{i}^{2}(\boldsymbol{%
\alpha },\sigma ^{2})\}_{i=1}^{n}).
\end{align*}%
For the heteroscedastic linear model, (\ref{eq1}) or (\ref{eq2}), the full
parameter vector is given by 
\begin{equation}
\boldsymbol{\theta }=(\boldsymbol{\alpha }^{T},\sigma ^{2},\boldsymbol{\beta 
}^{T})^{T},  \label{eq5}
\end{equation}%
and the hypothesis of homoscedasticity against heteroscedasticity
establishes 
\begin{equation}
H_{0}:\boldsymbol{\alpha }=\boldsymbol{0}_{r}\quad vs.\quad H_{1}:%
\boldsymbol{\alpha }\neq \boldsymbol{0}_{r}.  \label{eq5b}
\end{equation}

\begin{proposition}
\label{PropEsti}For the heteroscedastic linear model (\ref{eq1}) or (\ref%
{eq2}), the estimating function for the unrestricted minimum DPD estimators
of $\boldsymbol{\theta }$ is given by $\boldsymbol{U}_{n,\beta }(%
\boldsymbol{Y};\boldsymbol{\theta })=\boldsymbol{0}_{p+r+2}$, where%
\begin{equation*}
\boldsymbol{U}_{n,\beta }(\boldsymbol{\theta })=\frac{1}{n}\sum_{i=1}^{n}%
\boldsymbol{u}_{i,\beta }(Y_{i};\boldsymbol{\theta })=\frac{1}{n}%
\sum_{i=1}^{n}\left( \boldsymbol{u}_{i,\beta ,\boldsymbol{\alpha }%
}^{T}(Y_{i};\boldsymbol{\theta }),\boldsymbol{u}_{i,\beta ,\sigma
^{2}}^{T}(Y_{i};\boldsymbol{\theta }),\boldsymbol{u}_{i,\beta ,\boldsymbol{%
\beta }}^{T}(Y_{i};\boldsymbol{\theta })\right) ^{T},
\end{equation*}%
with%
\begin{align}
\boldsymbol{u}_{i,\beta ,\boldsymbol{\alpha }}(Y_{i};\boldsymbol{\theta })& =%
\frac{1}{(2\pi )^{\frac{\beta }{2}}(\sigma _{i}^{2}(\boldsymbol{\alpha }%
,\sigma ^{2}))^{\frac{\beta }{2}}}\frac{1}{2}\frac{h^{\prime }(\boldsymbol{z}%
_{i}^{T}\boldsymbol{\alpha })}{h(\boldsymbol{z}_{i}^{T}\boldsymbol{\alpha })}%
\left[ \exp \left\{ -\tfrac{\beta }{2}g_{i}(\boldsymbol{\theta })\right\}
\left( g_{i}(\boldsymbol{\theta })-1\right) +\frac{\beta }{(\beta +1)^{\frac{%
3}{2}}}\right] \boldsymbol{z}_{i},  \label{eq3} \\
u_{i,\beta ,\sigma ^{2}}(Y_{i};\boldsymbol{\theta })& =\frac{1}{(2\pi )^{%
\frac{\beta }{2}}(\sigma _{i}^{2}(\boldsymbol{\alpha },\sigma ^{2}))^{\frac{%
\beta }{2}}}\frac{1}{2\sigma ^{2}}\left[ \exp \left\{ -\tfrac{\beta }{2}%
g_{i}(\boldsymbol{\theta })\right\} \left( g_{i}(\boldsymbol{\theta }%
)-1\right) +\frac{\beta }{(\beta +1)^{\frac{3}{2}}}\right] ,  \notag \\
\boldsymbol{u}_{i,\beta ,\boldsymbol{\beta }}(Y_{i};\boldsymbol{\theta })& =%
\frac{1}{(2\pi )^{\frac{\beta }{2}}(\sigma _{i}^{2}(\boldsymbol{\alpha }%
,\sigma ^{2}))^{\frac{\beta }{2}}}\exp \left\{ -\tfrac{\beta }{2}g_{i}(%
\boldsymbol{\theta })\right\} \frac{1}{\sigma ^{2}}\frac{\epsilon _{i}(%
\boldsymbol{\beta })}{h(\boldsymbol{z}_{i}^{T}\boldsymbol{\alpha })}%
\boldsymbol{x}_{i},  \notag
\end{align}%
with%
\begin{align}
	g_{i}(\boldsymbol{\theta}) &  =\frac{\epsilon_{i}^{2}(\boldsymbol{\beta}%
		)}{\sigma_{i}^{2}(\boldsymbol{\alpha},\sigma^{2})},\label{eq6}\\
	\epsilon_{i} &  =\epsilon_{i}(\boldsymbol{\beta})=Y_{i}-\boldsymbol{x}_{i}%
	^{T}\boldsymbol{\beta}.\label{error}%
\end{align}
\end{proposition}

\begin{proof}
From Mart\'{\i}n (2021), taking $p=1$ and Theorem 6, the parameter vector $%
\boldsymbol{\eta }_{i}=\boldsymbol{\eta }_{i}(\boldsymbol{\theta })=(\sigma
_{i}^{2},\mu _{i})^{T}$, with $\mu _{i}=\mu _{i}(\boldsymbol{\beta })=%
\boldsymbol{x}_{i}^{T}\boldsymbol{\beta }$, $\sigma _{i}^{2}=\sigma ^{2}h(%
\boldsymbol{z}_{i}^{T}\boldsymbol{\alpha })$ is taken into account. By
following the chain rule of differentiation, we get%
\begin{equation*}
\boldsymbol{u}_{i}(Y_{i};\boldsymbol{\theta })=\frac{\partial }{\partial 
\boldsymbol{\theta }}\boldsymbol{\eta }_{i}^{T}(\boldsymbol{\theta })%
\boldsymbol{u}_{i}(Y_{i};\boldsymbol{\eta }_{i,\beta }),
\end{equation*}%
where $\boldsymbol{u}_{i}(Y_{i};\boldsymbol{\eta }_{i,\beta
})=(u_{i}(Y_{i};\sigma _{i}^{2}),u_{i}(Y_{i};\mu _{i}))^{T}$,%
\begin{equation}
\frac{\partial }{\partial \boldsymbol{\theta }}\boldsymbol{\eta }_{i}^{T}(%
\boldsymbol{\theta })=%
\begin{pmatrix}
\sigma ^{2}h^{\prime }(\boldsymbol{z}_{i}^{T}\boldsymbol{\alpha })%
\boldsymbol{z}_{i} & 0 \\ 
h(\boldsymbol{z}_{i}^{T}\boldsymbol{\alpha }) & 0 \\ 
0 & \boldsymbol{x}_{i}%
\end{pmatrix}%
,  \label{eq7}
\end{equation}%
and finally 
\begin{align*}
\boldsymbol{u}_{i,\beta ,\boldsymbol{\alpha }}(Y_{i};\boldsymbol{\theta })&
=\sigma ^{2}h^{\prime }(\boldsymbol{z}_{i}^{T}\boldsymbol{\alpha }%
)u_{i}(Y_{i};\sigma _{i}^{2})\boldsymbol{z}_{i}, \\
u_{i,\beta ,\sigma ^{2}}(Y_{i};\boldsymbol{\theta })& =h(\boldsymbol{z}%
_{i}^{T}\boldsymbol{\alpha })u_{i}(Y_{i};\sigma _{i}^{2}), \\
\boldsymbol{u}_{i,\beta ,\boldsymbol{\beta }}(Y_{i};\boldsymbol{\theta })&
=u_{i}(Y_{i};\mu _{i})\boldsymbol{x}_{i},
\end{align*}%
with%
\begin{align*}
& \\
u_{i}(Y_{i};\sigma _{i}^{2})& =\frac{1}{(2\pi )^{\frac{\beta }{2}}(\sigma
_{i}^{2}(\boldsymbol{\alpha },\sigma ^{2}))^{\frac{\beta }{2}}}\frac{1}{%
2\sigma ^{2}h(\boldsymbol{z}_{i}^{T}\boldsymbol{\alpha })}\left[ \exp
\left\{ -\tfrac{\beta }{2}g_{i}(\boldsymbol{\theta })\right\} \left( g_{i}(%
\boldsymbol{\theta })-1\right) +\frac{\beta }{(\beta +1)^{\frac{3}{2}}}%
\right] , \\
u_{i}(Y_{i};\mu _{i})& =\frac{1}{(2\pi )^{\frac{\beta }{2}}(\sigma _{i}^{2}(%
\boldsymbol{\alpha },\sigma ^{2}))^{\frac{\beta }{2}}}\exp \left\{ -\tfrac{%
\beta }{2}g_{i}(\boldsymbol{\theta })\right\} \frac{1}{\sigma ^{2}}\frac{%
\epsilon _{i}(\boldsymbol{\beta })}{h(\boldsymbol{z}_{i}^{T}\boldsymbol{%
\alpha })}.
\end{align*}
\end{proof}

\begin{corollary}
	\label{Corollary2}For the heteroscedastic linear model, (\ref{eq1}) or
	(\ref{eq2}), the restricted minimum DPD estimators of $\boldsymbol{\theta}$,
	under homoscedastic null hypothesis, $\widetilde{\boldsymbol{\theta}}_{\beta
	}=(\boldsymbol{0}_{r}^{T},\widetilde{\sigma}_{\beta}^{2}%
	,\widetilde{\boldsymbol{\beta}}_{\beta}^{T})^{T}$, is obtained as solution of
	\begin{align}
		\mathbb{X}_{n}^{T}\mathrm{diag}\left(  \exp(-\tfrac{\beta}{2}\boldsymbol{g}%
		(\widetilde{\boldsymbol{\theta}}_{\beta}))\right)  \left(  \boldsymbol{Y}%
		-\mathbb{X}_{n}\widetilde{\boldsymbol{\beta}}_{\beta}\right)   &
		=\boldsymbol{0}_{r},\label{eqAlgo2}\\
		\boldsymbol{1}_{n}^{T}\boldsymbol{v}(\widetilde{\boldsymbol{\theta}}_{\beta})
		&  =0,\label{eqAlgo1}%
	\end{align}
	where
	\begin{align}
		\boldsymbol{g}(\widetilde{\boldsymbol{\theta}}_{\beta}) &  =(g_{1}%
		(\widetilde{\boldsymbol{\theta}}_{\beta}),\ldots,g_{n}%
		(\widetilde{\boldsymbol{\theta}}_{\beta}))^{T},\nonumber\\
		g_{i}(\widetilde{\boldsymbol{\theta}}_{\beta}) &  =\frac{\epsilon_{i}%
			^{2}(\widetilde{\boldsymbol{\beta}}_{\beta})}{\widetilde{\sigma}_{\beta}^{2}%
		},\quad\epsilon_{i}(\widetilde{\boldsymbol{\beta}}_{\beta})=Y_{i}%
		-\boldsymbol{x}_{i}^{T}\widetilde{\boldsymbol{\beta}}_{\beta},\quad
		i=1,\ldots,n,\nonumber\\
		\exp(-\tfrac{\beta}{2}\boldsymbol{g}(\widetilde{\boldsymbol{\theta}}_{\beta}))
		&  =(\exp\{-\tfrac{\beta}{2}g_{1}(\widetilde{\boldsymbol{\theta}}_{\beta
		})\},\ldots,\exp\{-\tfrac{\beta}{2}g_{n}(\widetilde{\boldsymbol{\theta}%
		}_{\beta})\})^{T},\label{eqAlgo3}\\
		\boldsymbol{v}(\widetilde{\boldsymbol{\theta}}_{\beta}) &  =\mathrm{diag}%
		\left(  \exp(-\tfrac{\beta}{2}\boldsymbol{g}(\widetilde{\boldsymbol{\theta}%
		}_{\beta}))\right)  (\boldsymbol{g}(\widetilde{\boldsymbol{\theta}}_{\beta
		})-\boldsymbol{1}_{n})+\frac{\beta}{(\beta+1)^{\frac{3}{2}}}\boldsymbol{1}%
		_{n}.\label{vectorLM}%
	\end{align}
	In addition, the corresponding restricted minimum DPD estimators of the LM
	vector is given by
	\begin{equation}
		\boldsymbol{\lambda}_{n,\beta}(\widetilde{\boldsymbol{\theta}}_{\beta}%
		)=\frac{1}{(2\pi)^{\frac{\beta}{2}}(\widetilde{\sigma}_{\beta}^{2}%
			)^{\frac{\beta}{2}}}\frac{h^{\prime}(0)}{2}\frac{1}{n}\mathbb{Z}_{n}%
		^{T}\boldsymbol{v}(\widetilde{\boldsymbol{\theta}}_{\beta}),\label{eq3b}%
	\end{equation}
	with
	\[
	\mathbb{Z}_{n}=(\boldsymbol{z}_{1},\ldots,\boldsymbol{z}_{n})^{T},
	\]
	having $r$ columns ($n$-vectors) linearly independent of $\boldsymbol{1}_{n}$.
\end{corollary}

\begin{proof}
	(\ref{eqAlgo2}) is a direct from Proposition \ref{PropEsti}, taking
	$\boldsymbol{\alpha=0}_{r}$ and $\boldsymbol{u}_{i,\beta,\boldsymbol{\beta}%
	}(Y_{i};\boldsymbol{\theta})=\boldsymbol{0}_{p+1}$, while (\ref{eq3b}) comes
	from
	\[
	\sum\limits_{i=1}^{n}\boldsymbol{u}_{i,\beta,\boldsymbol{\alpha=0}_{p}}%
	(Y_{i};\widetilde{\boldsymbol{\theta}}_{\beta})=\frac{\frac{h^{\prime}(0)}{2}%
	}{(2\pi)^{\frac{\beta}{2}}(\widetilde{\sigma}_{\beta}^{2})^{\frac{\beta}{2}}%
	}\sum\limits_{i=1}^{n}\left(  \exp\left\{  -\tfrac{\beta}{2}g_{i}%
	(\widetilde{\boldsymbol{\theta}}_{\beta})\right\}  \left(  g_{i}%
	(\widetilde{\boldsymbol{\theta}}_{\beta})-1\right)  +\frac{\beta}%
	{(\beta+1)^{\frac{3}{2}}}\right)  \boldsymbol{z}_{i}%
	\]
	and $\boldsymbol{\lambda}_{n,\beta}(\widetilde{\boldsymbol{\theta}}_{\beta
	})=\frac{1}{n}\sum\limits_{i=1}^{n}\boldsymbol{u}_{i,\beta,\boldsymbol{\alpha
			=0}_{p}}(Y_{i};\widetilde{\boldsymbol{\theta}}_{\beta})$ , according to
	(\ref{eqA}).
\end{proof}

\section{Main theoretical results\label{Sec2}}

In this section we shall focus in a general model such as the one introduced in Section \ref{Sec1}. For the particular case of identically distributed observations, the
asymptotic distribution of $\sqrt{n}(\widetilde{\boldsymbol{\theta }}_{\beta
}-\boldsymbol{\theta }_{0})$ was known, from Basu et al. (2017). The
following result generalizes the previous one in two ways, extending to
non-identically distributed observations, and also by considering jointly
the estimator of the Lagrange multipliers. In what is to follow, it is assumed to be fulfilled
the regularity conditions given in Basu et al. (2017) based on a single
observation, as well as the existence of the limits of the following two
matrices based on the whole set of observations,%
\begin{equation}
\boldsymbol{J}_{\beta }(\boldsymbol{\theta })=\lim_{n\rightarrow \infty }%
\boldsymbol{\bar{J}}_{n,\beta }(\boldsymbol{\theta })=\lim_{n\rightarrow
\infty }\frac{1}{n}\sum_{i=1}^{n}\boldsymbol{J}_{i,\beta }(\boldsymbol{%
\theta }),  \label{Jmatrix}
\end{equation}%
with 
\begin{equation*}
\boldsymbol{J}_{i,\beta }(\boldsymbol{\theta })=-\mathrm{E}\left[ \frac{%
\partial }{\partial \boldsymbol{\theta }}\boldsymbol{u}_{i,\beta }^{T}\left(
Y_{i},\boldsymbol{\theta }\right) \right] =\int_{\mathcal{-\infty }%
}^{+\infty }f_{i,\boldsymbol{\theta }}^{\beta +1}(y)\boldsymbol{s}_{i,%
\boldsymbol{\theta }}(y)\boldsymbol{s}_{i,\boldsymbol{\theta }}^{T}(y)dy=%
\mathrm{E}[f_{i,\boldsymbol{\theta }}^{\beta }(Y_{i})\boldsymbol{s}_{i,%
\boldsymbol{\theta }}(Y_{i})\boldsymbol{s}_{i,\boldsymbol{\theta }%
}^{T}(Y_{i})],
\end{equation*}%
and%
\begin{equation}
\boldsymbol{K}_{\beta }\left( \boldsymbol{\theta }\right)
=\lim_{n\rightarrow \infty }\boldsymbol{\bar{K}}_{n,\beta }(\boldsymbol{%
\theta })=\lim_{n\rightarrow \infty }\frac{1}{n}\sum_{i=1}^{n}\boldsymbol{K}%
_{i,\beta }\left( \boldsymbol{\theta }\right) ,  \label{Kmatrix}
\end{equation}%
with 
\begin{align*}
\boldsymbol{\bar{K}}_{n,\beta }(\boldsymbol{\theta })& =\frac{1}{n}%
\sum_{i=1}^{n}\boldsymbol{K}_{i,\beta }\left( \boldsymbol{\theta }\right) =%
\boldsymbol{\bar{J}}_{n,2\beta }(\boldsymbol{\theta })-\frac{1}{n}%
\sum_{i=1}^{n}\boldsymbol{\xi }_{i,\beta }\left( \boldsymbol{\theta }\right) 
\boldsymbol{\xi }_{i,\beta }^{T}\left( \boldsymbol{\theta }\right) , \\
\boldsymbol{K}_{i,\beta }(\boldsymbol{\theta })& =\mathrm{E}\left[ 
\boldsymbol{u}_{i,\beta }\left( Y_{i},\boldsymbol{\theta }\right) 
\boldsymbol{u}_{i,\beta }^{T}\left( Y_{i},\boldsymbol{\theta }\right) \right]
=\boldsymbol{J}_{i,2\beta }(\boldsymbol{\theta })-\boldsymbol{\xi }_{i,\beta
}(\boldsymbol{\theta })\boldsymbol{\xi }_{i,\beta }^{T}(\boldsymbol{\theta }%
),
\end{align*}%
and $\boldsymbol{\xi }_{i,\beta }(\boldsymbol{\theta })$ given by (\ref{psi}%
).

\begin{theorem}
\label{ThM}Restricted to the null hypothesis established by (\ref{eq4}),
under the existence assumption of (\ref{Jmatrix})-(\ref{Kmatrix}), it holds%
\begin{equation*}
\sqrt{n}%
\begin{bmatrix}
\widetilde{\boldsymbol{\theta }}_{\beta }-\boldsymbol{\theta }_{0} \\ 
\boldsymbol{\lambda }_{n,\beta }(\widetilde{\boldsymbol{\theta }}_{\beta })%
\end{bmatrix}%
\underset{n\rightarrow \infty }{\overset{\mathcal{L}}{\longrightarrow }}%
\mathcal{N}\left( \boldsymbol{0}_{p},%
\begin{bmatrix}
\boldsymbol{P}_{\beta }(\boldsymbol{\theta }_{0}) & \boldsymbol{R}_{\beta }(%
\boldsymbol{\theta }_{0}) \\ 
\boldsymbol{R}_{\beta }^{T}(\boldsymbol{\theta }_{0}) & \boldsymbol{Q}%
_{\beta }(\boldsymbol{\theta }_{0})%
\end{bmatrix}%
\right) ,
\end{equation*}%
where the normal distribution has singular variance-covariance with%
\begin{align}
\boldsymbol{P}_{\beta }(\boldsymbol{\theta })& =\boldsymbol{\Sigma }%
_{11,\beta }(\boldsymbol{\theta })\boldsymbol{K}_{\beta }(\boldsymbol{\theta 
})\boldsymbol{\Sigma }_{11,\beta }^{T}(\boldsymbol{\theta }),  \label{v1} \\
\boldsymbol{Q}_{\beta }(\boldsymbol{\theta })& =\boldsymbol{\Sigma }%
_{21,\beta }(\boldsymbol{\theta })\boldsymbol{K}_{\beta }(\boldsymbol{\theta 
})\boldsymbol{\Sigma }_{21,\beta }^{T}(\boldsymbol{\theta }),  \label{v2} \\
\boldsymbol{R}_{\beta }(\boldsymbol{\theta })& =\boldsymbol{\Sigma }%
_{11,\beta }(\boldsymbol{\theta })\boldsymbol{K}_{\beta }(\boldsymbol{\theta 
})\boldsymbol{\Sigma }_{21,\beta }^{T}(\boldsymbol{\theta }),  \label{v3}
\end{align}%
and%
\begin{align*}
\boldsymbol{\Sigma }_{11,\beta }(\boldsymbol{\theta })& =\boldsymbol{J}%
_{\beta }^{-1}(\boldsymbol{\theta })-\boldsymbol{J}_{\beta }^{-1}(%
\boldsymbol{\theta })\boldsymbol{M}(\boldsymbol{\theta })\boldsymbol{\Sigma }%
_{21,\beta }(\boldsymbol{\theta }), \\
\boldsymbol{\Sigma }_{21,\beta }(\boldsymbol{\theta })& =\left[ \boldsymbol{M%
}^{T}(\boldsymbol{\theta })\boldsymbol{J}_{\beta }^{-1}(\boldsymbol{\theta })%
\boldsymbol{M}(\boldsymbol{\theta })\right] ^{-1}\boldsymbol{M}^{T}(%
\boldsymbol{\theta })\boldsymbol{J}_{\beta }^{-1}(\boldsymbol{\theta }), \\
.\boldsymbol{\Sigma }_{22,\beta }(\boldsymbol{\theta })& =\left[ \boldsymbol{%
M}^{T}(\boldsymbol{\theta })\boldsymbol{J}_{\beta }^{-1}(\boldsymbol{\theta }%
)\boldsymbol{M}(\boldsymbol{\theta })\right] ^{-1}.
\end{align*}
\end{theorem}

Appendix \ref{AppA} covers the proof of Theorem \ref{ThM}. The
original and classical Lagrange multipliers test, published by Aitchison and
Silvey (1958) and Silvey (1959), was for identically distributed
observation, in which taking into account that for $\beta =0$, $\boldsymbol{%
\bar{K}}_{n,\beta =0}(\boldsymbol{\theta })=\boldsymbol{\bar{J}}_{n,\beta
=0}(\boldsymbol{\theta })=\boldsymbol{\bar{I}}_{F,n}(\boldsymbol{\theta })$
is the average Information matrix, the role of these matrices was played by
the Information matrix based on a unique observation. The asymptotic
distribution of $\sqrt{n}\boldsymbol{\lambda }_{n,\beta }(\widetilde{%
\boldsymbol{\theta }}_{\beta })$ is the cornerstone for the following
definition and generalizes the one given in Basu et al. (2022), from which
it is not possible to derive the Breusch-Pagan $\beta $-score LM tests
presented in Section \ref{Sec3} of the current paper.

\begin{definition}
\label{def}The $\beta $-score Lagrange multipliers test for non-identically
distributed individuals is given by%
\begin{equation}
R_{n,\beta }(\widetilde{\boldsymbol{\theta }}_{\beta })=n\boldsymbol{U}%
_{n,\beta }^{T}(\widetilde{\boldsymbol{\theta }}_{\beta })\boldsymbol{\bar{J}%
}_{n,\beta }^{-1}(\widetilde{\boldsymbol{\theta }}_{\beta })\boldsymbol{M}(%
\widetilde{\boldsymbol{\theta }}_{\beta })\left[ \boldsymbol{M}^{T}(%
\widetilde{\boldsymbol{\theta }}_{\beta })\boldsymbol{\bar{V}}_{n,\beta }(%
\widetilde{\boldsymbol{\theta }}_{\beta })\boldsymbol{M}(\widetilde{%
\boldsymbol{\theta }}_{\beta })\right] ^{-1}\boldsymbol{M}^{T}(\widetilde{%
\boldsymbol{\theta }}_{\beta })\boldsymbol{\bar{J}}_{n,\beta }^{-1}(%
\widetilde{\boldsymbol{\theta }}_{\beta })\boldsymbol{U}_{n,\beta }(%
\widetilde{\boldsymbol{\theta }}_{\beta }),  \label{LM_TS}
\end{equation}%
where%
\begin{equation}
\boldsymbol{\bar{V}}_{n,\beta }(\boldsymbol{\theta })=\boldsymbol{\bar{J}}%
_{n,\beta }^{-1}(\boldsymbol{\theta })\boldsymbol{\bar{K}}_{n,\beta }(%
\boldsymbol{\theta })\boldsymbol{\bar{J}}_{n,\beta }^{-1}(\boldsymbol{\theta 
})  \label{Qbar}
\end{equation}%
is the \textquotedblleft empirical sandwich matrix\textquotedblright .
\end{definition}

\begin{theorem}
\label{ThM2}The asymptotic distribution of (\ref{LM_TS}) is a chi-square
with $r$ degrees of freedom.
\end{theorem}

Large values of the test statistic, on the right hand side tail of $\chi
_{r}^{2}$, are interpreted as strong evidence against the null hypothesis.
Appendix \ref{AppB} covers the proof of Theorem \ref{ThM2}.

\begin{theorem}
	\label{contiguous}Consider the test $H_{0}:\boldsymbol{\theta}\in\Theta_{0}%
	$\ vs.$\ H_{1}:\boldsymbol{\theta}\in\Theta-\Theta_{0}$, with $\Theta_{0}$
	given by (\ref{eq3}) and let a sequence of local Pitman-type alternatives be
	defined by $H_{1,n}:\boldsymbol{\theta}_{n}\in\Theta_{1,n}$, with
	$\Theta_{1,n}$ given by
	\begin{equation}
		\Theta_{1,n}=\left\{  \boldsymbol{\theta}_{n}\in\Theta:\boldsymbol{m}%
		(\boldsymbol{\theta})=\tfrac{1}{\sqrt{n}}\boldsymbol{\delta}\right\}
		,\label{pitman}%
	\end{equation}
	being fixed $\boldsymbol{\delta}\in%
	\mathbb{R}
	^{r}-\{\boldsymbol{0}_{r}\}$. Under the sequence $H_{1,n}:\boldsymbol{\theta
	}\in\Theta_{1,n}$, the asymptotic distribution of the $\beta$-score LM
	test-statistic for non-identically distributed individuals, (\ref{LM_TS}), is
	a chi-square with $r$ degrees of freedom and non-centrality parameter
	\begin{equation}
		\nu_{\beta}(\boldsymbol{\theta}_{0},\boldsymbol{\delta})=\boldsymbol{\delta
		}^{T}\left[  \boldsymbol{M}^{T}(\boldsymbol{\theta}_{0})\boldsymbol{V}_{\beta
		}(\boldsymbol{\theta}_{0})\boldsymbol{M}(\boldsymbol{\theta}_{0})\right]
		^{-1}\boldsymbol{\delta},\label{NCP}%
	\end{equation}
	where
	\begin{equation}
		\boldsymbol{V}_{\beta}(\boldsymbol{\theta})=\boldsymbol{J}_{\beta}%
		^{-1}(\boldsymbol{\theta})\boldsymbol{K}_{\beta}(\boldsymbol{\theta
		})\boldsymbol{J}_{\beta}^{-1}(\boldsymbol{\theta})\label{v}%
	\end{equation}
	is the \textquotedblleft theoretical sandwich matrix\textquotedblright, the
	asymptotic variance-covariance matrix of $\widetilde{\boldsymbol{\theta}%
	}_{\beta}$. The corresponding asymptotic power function, given the nominal
	level $\alpha$, is
	\begin{equation}
		\pi_{\beta}(\boldsymbol{y};\boldsymbol{\theta}_{0},\boldsymbol{\delta}%
		)=\lim_{n\rightarrow\infty}P(R_{n,\beta}(\widetilde{\boldsymbol{\theta}%
		}_{\beta})>\chi_{r,\alpha}^{2}|\underline{F}_{n}(\boldsymbol{y}%
		;\boldsymbol{\theta}_{0},\boldsymbol{\delta}))=Q_{\frac{r}{2}}(\nu_{\beta
		}(\boldsymbol{\theta}_{0},\boldsymbol{\delta}),\chi_{r,\alpha}^{2}%
		),\label{power}%
	\end{equation}
	where $Q_{M}(a,b)$ is the Marcum $Q$-function (see details in Nuttal, 1975).
\end{theorem}

From the previous result, the Breusch-Pagan test is consistent under Pitman
alternatives since taking $\boldsymbol{\delta}^{\ast}=\sqrt{n}%
\boldsymbol{\delta}$, we get $\lim_{\nu\rightarrow\infty}Q_{\frac{r}{2}}%
(\nu,\chi_{r,\alpha}^{2})=1$.

The Pitman asymptotic relative efficiency (ARE) of the $\beta $-score LM
test to the classical LM test is given by 
\begin{equation}
\mathrm{ARE}(R_{n,\beta }(\widetilde{\boldsymbol{\theta }}_{\beta
}),R_{n,\beta =0}(\widetilde{\boldsymbol{\theta }}_{\beta =0});\boldsymbol{%
\delta })=\frac{\nu _{\beta =0}(\boldsymbol{\theta }_{0},\boldsymbol{\delta }%
)}{\nu _{\beta }(\boldsymbol{\theta }_{0},\boldsymbol{\delta })}>1,
\label{AREg}
\end{equation}%
based on similar ideas taken from Hannan (1956) and Koenker and Bassett
(1982).

As particular case of Definition \ref{def}, we get the following version.

\begin{proposition}
	\label{ThM2b}Let us consider the particular case of $\boldsymbol{\theta
	}=(\boldsymbol{\theta}_{1}^{T},\boldsymbol{\theta}_{2}^{T},\boldsymbol{\theta
	}_{3}^{T})^{T}$ and $\boldsymbol{m}(\boldsymbol{\theta})=\boldsymbol{\theta
	}_{1}-\boldsymbol{\theta}_{1,0}$, with $(\widetilde{\boldsymbol{\theta}%
	}_{1,\beta}^{T},\widetilde{\boldsymbol{\theta}}_{2,\beta}^{T})^{T}$ and
	$\widetilde{\boldsymbol{\theta}}_{3,\beta}$ being asymptotically independent
	and $\boldsymbol{\bar{K}}_{n,\beta}(\boldsymbol{\theta})$ is a proportional
	matrix, with respect to $\boldsymbol{\bar{J}}_{n,\beta}(\boldsymbol{\theta})$,
	within the block correspondent to $(\widetilde{\boldsymbol{\theta}}_{1,\beta
	}^{T},\widetilde{\boldsymbol{\theta}}_{2,\beta}^{T})^{T}$. Then, the $\beta
	$-score Lagrange multipliers test for non-identically distributed individuals,
	given in (\ref{LM_TS}), has the following simpler expression
	\begin{equation}
		R_{n,\beta}(\widetilde{\boldsymbol{\theta}}_{\beta})=n\boldsymbol{\lambda
		}_{n,\beta}^{T}(\widetilde{\boldsymbol{\theta}}_{\beta})\left(
		\boldsymbol{\bar{K}}_{n,\beta}^{-1}(\widetilde{\boldsymbol{\theta}}_{\beta
		})\right)  _{11}\boldsymbol{\lambda}_{n,\beta}(\widetilde{\boldsymbol{\theta}%
		}_{\beta}).\label{LM_TS2}%
	\end{equation}
	where
	\[
	\left(  \boldsymbol{\bar{K}}_{n,\beta}^{-1}(\widetilde{\boldsymbol{\theta}%
	}_{\beta})\right)  _{11}=(\boldsymbol{I}_{r},\boldsymbol{0}_{r\times
		(p-r)})\boldsymbol{\bar{K}}_{n,\beta}^{-1}(\boldsymbol{\theta})(\boldsymbol{I}%
	_{r},\boldsymbol{0}_{r\times(p-r)})^{T}.
	\]
	
\end{proposition}

\section{Derivation of the Breusch-Pagan $\protect\beta $-score LM tests 
\label{Sec3}}

Focussed on the matrix version of Definition \ref{def}, let us consider the extended design matrix of the scedasticity function,
$\mathbb{\breve{Z}}_{n}=(1_{n},\mathbb{Z}_{n})$ and the corresponding rows,
$\boldsymbol{\breve{z}}_{i}=(1,\boldsymbol{z}_{i})^{T}$, $i=1,\ldots,n$. We shall assume that
$\Lambda_{\min}(\frac{1}{n}\mathbb{X}_{n}^{T}\mathbb{X}_{n})>0$,
$\Lambda_{\min}(\frac{1}{n}\mathbb{\breve{Z}}_{n}^{T}\mathbb{\breve{Z}}%
_{n})>0$ and%
\begin{equation}
	\lim_{n\rightarrow\infty}\frac{1}{\sqrt{n}}\frac{\underset{1\leq i\leq
			n}{\max}\left\Vert \boldsymbol{x}_{i}\right\Vert }{\Lambda_{\min}%
		^{\frac{1}{2}}(\frac{1}{n}\mathbb{X}_{n}^{T}\mathbb{X}_{n})}=\lim
	_{n\rightarrow\infty}\frac{1}{\sqrt{n}}\frac{\underset{1\leq i\leq n}{\max
		}\left\Vert \boldsymbol{\breve{z}}_{i}\right\Vert}{\Lambda_{\min}%
		^{\frac{1}{2}}(\frac{1}{n}\mathbb{\breve{Z}}_{n}^{T}\mathbb{\breve{Z}}_{n}%
		)}=0,\label{cond}%
\end{equation}
where $\Lambda_{\min}(\cdot)$ denotes the minimum eigenvalue of a matrix.

Using the full parameter given in (\ref{eq5}), the hypothesis of
homoscedasticity against heteroscedasticity, (\ref{eq5b}), belongs to the
particular case of Theorem \ref{ThM2b} with $\boldsymbol{\alpha }$ playing
the role of $\boldsymbol{\theta }_{1}$ and $(\sigma ^{2},\boldsymbol{\beta }%
^{T})^{T}$ the role of $\boldsymbol{\theta }_2$. The current section is
mainly devoted to calculate $\boldsymbol{\bar{K}}_{n,\beta ,\boldsymbol{%
\alpha }}(\widetilde{\boldsymbol{\theta }}_{\beta })$ and $\boldsymbol{\bar{J%
}}_{n,\beta ,\boldsymbol{\alpha }}^{-1}(\widetilde{\boldsymbol{\theta }}%
_{\beta })$.

\begin{theorem}
\label{ThM3}For the heteroscedastic linear model, (\ref{eq1}) or (\ref{eq2}),
under homoscedastic null hypothesis, it holds%
\begin{align}
	\boldsymbol{\bar{K}}_{n,\beta }(\widetilde{\boldsymbol{\theta }}_{\beta })& =%
	\begin{pmatrix}
		a_{1,\beta }\boldsymbol{W}_{n,\beta }(\widetilde{\boldsymbol{\theta }}%
		_{\beta }) & \boldsymbol{0}_{(r+1)\times (p+1)} \\ 
		\boldsymbol{0}_{(p+1)\times (r+1)} & a_{2,\beta }\frac{1}{n}\mathbb{X}%
		_{n}^{T}\mathbb{X}_{n}%
	\end{pmatrix}%
	,  \label{K} \\
	\boldsymbol{\bar{J}}_{n,\beta }(\widetilde{\boldsymbol{\theta }}_{\beta })& =%
	\begin{pmatrix}
		b_{1,\beta }\boldsymbol{W}_{n,\beta }(\widetilde{\boldsymbol{\theta }}%
		_{\beta }) & \boldsymbol{0}_{(r+1)\times (p+1)} \\ 
		\boldsymbol{0}_{(p+1)\times (r+1)} & b_{2,\beta }\frac{1}{n}\mathbb{X}%
		_{n}^{T}\mathbb{X}_{n}%
	\end{pmatrix}%
	,  \label{J} \\
	\boldsymbol{W}_{n,\beta }(\widetilde{\boldsymbol{\theta }}_{\beta })& =%
	\begin{pmatrix}
		\widetilde{\sigma }_{\beta }^{4}(h^{\prime }(0))^{2}\frac{1}{n}\mathbb{Z}%
		_{n}^{T}\mathbb{Z}_{n} & \widetilde{\sigma }_{\beta }^{2}h^{\prime }(0)%
		\overline{\boldsymbol{z}}_{n} \\ 
		\widetilde{\sigma }_{\beta }^{2}h^{\prime }(0)\overline{\boldsymbol{z}}%
		_{n}^{T} & 1%
	\end{pmatrix}%
	,  \label{W}
\end{align}%
where $\frac{1}{n}\mathbb{X}_{n}^{T}\mathbb{X}_{n}\mathbb{=}\frac{1}{n}%
\sum_{i=1}^{n}\boldsymbol{x}_{i}\boldsymbol{x}_{i}^{T}$, $\frac{1}{n}\mathbb{%
Z}_{n}^{T}\mathbb{Z}_{n}\mathbb{=}\frac{1}{n}\sum_{i=1}^{n}\boldsymbol{z}_{i}%
\boldsymbol{z}_{i}^{T}$, $\overline{\boldsymbol{z}}_{n}=\frac{1}{n}\mathbb{Z}%
_{n}^{T}\boldsymbol{1}_{n}$,
\begin{align*}
	a_{1,\beta }& =\frac{\frac{2\beta ^{2}+1}{2(2\beta +1)^{5/2}}-\frac{\beta
			^{2}}{4(\beta +1)^{3}}}{(2\pi )^{\beta }(\widetilde{\sigma }_{\beta
		}^{2})^{\beta +2}},\quad a_{2,\beta }=\frac{1}{(2\pi )^{\beta }(\widetilde{%
			\sigma }_{\beta }^{2})^{\beta +2}(2\beta +1)^{\frac{3}{2}}}, \\
	b_{1,\beta }& =\frac{\frac{\beta ^{2}+2}{4(\beta +1)^{5/2}}}{(2\pi )^{\frac{%
				\beta }{2}}(\widetilde{\sigma }_{\beta }^{2})^{\frac{\beta }{2}+2}},\quad
	b_{2,\beta }=\frac{1}{(2\pi )^{\frac{\beta }{2}}(\widetilde{\sigma }_{\beta
		}^{2})^{\frac{\beta }{2}+1}(\beta +1)^{\frac{3}{2}}}.
\end{align*}
\end{theorem}

\begin{proof}
From (\ref{eq7}), taking $h(0)=1$, we get%
\begin{align}
\boldsymbol{K}_{i,\beta }(\widetilde{\boldsymbol{\theta }}_{\beta })& =%
\mathrm{E}[\boldsymbol{u}_{i,\boldsymbol{\theta }}(Y_{i};\widetilde{%
\boldsymbol{\theta }}_{\beta })\boldsymbol{u}_{i,\boldsymbol{\theta }%
}^{T}(Y_{i};\widetilde{\boldsymbol{\theta }}_{\beta })]  \notag \\
& =\left. \frac{\partial }{\partial \boldsymbol{\theta }}\boldsymbol{\eta }%
_{i}^{T}(\boldsymbol{\theta })\right] _{\boldsymbol{\theta }=\widetilde{%
\boldsymbol{\theta }}_{\beta }}\boldsymbol{K}_{\beta }(\boldsymbol{\eta }%
_{i}(\widetilde{\boldsymbol{\theta }}_{\beta }))\left. \frac{\partial }{%
\partial \boldsymbol{\theta }^{T}}\boldsymbol{\eta }_{i}(\boldsymbol{\theta }%
)\right] _{\boldsymbol{\theta }=\widetilde{\boldsymbol{\theta }}_{\beta }},
\label{eq9}
\end{align}%
where%
\begin{equation}
	\boldsymbol{K}_{\beta }(\boldsymbol{\eta }_{i}(\widetilde{\boldsymbol{\theta 
	}}_{\beta }))=\frac{1}{(2\pi )^{\beta }(\widetilde{\sigma }_{\beta
		}^{2})^{\beta +2}}%
	\begin{pmatrix}
		\frac{2\beta ^{2}+1}{2(2\beta +1)^{5/2}}-\frac{\beta ^{2}}{4(\beta +1)^{3}}
		& 0 \\ 
		0 & \frac{\widetilde{\sigma }_{\beta }^{2}}{(2\beta +1)^{\frac{3}{2}}}%
	\end{pmatrix}%
	,  \label{eq10}
\end{equation}%
can be obtained from Mart\'{\i}n 2020, Corollary 4. Finally, from $\boldsymbol{\bar{K}}_{\beta }(\widetilde{\boldsymbol{\theta }}%
_{\beta })=\frac{1}{n}\sum\limits_{i=1}^{n}\boldsymbol{K}_{i,\beta }(%
\widetilde{\boldsymbol{\theta }}_{\beta })$ the desired expression of (\ref%
{K}) is\ obtained. The expression of (\ref{J}) is\ obtained in a similar way.
\end{proof}
\begin{remark}
	The assumptions given in (\ref{cond}) arise as application of the multivariate
	Lindeberg Central Limit Theorem (CLT) to the estimators obtained as solution of the
	estimating equations given in Proposition \ref{PropEsti}, with weaker
	assumptions for the errors, just only considering null mean and finite variance.
\end{remark}

\begin{theorem}
\label{pitman2}The Pitman ARE of the Breusch-Pagan $\beta $-score LM
test-statistic (with respect to the classical one, $\beta =0$) for the
heteroscedastic linear model (\ref{eq1}) or (\ref{eq2}) is given by%
\begin{equation}
\mathrm{ARE}(R_{n,\beta }(\widetilde{\boldsymbol{\theta }}_{\beta
}),R_{n,\beta =0}(\widetilde{\boldsymbol{\theta }}_{\beta =0}))=\frac{8(\beta+1)^{5}}{(\beta^{2}+2)^{2}}\left(  \frac{2\beta^{2}+1}%
{2(2\beta+1)^{5/2}}-\frac{\beta^{2}}{4(\beta+1)^{3}}\right),\quad \beta >0.
\label{ARE}
\end{equation}
\end{theorem}

\begin{proof}
From Theorem \ref{ThM3} and equation (\ref{AREg}), we get the ARE to be $2b_{1,\beta}^{2}/a_{1,\beta}$.\medskip
\end{proof}

\begin{figure}[ptb]
	\centering
	\includegraphics[width=0.5\linewidth,height=0.45\textheight]{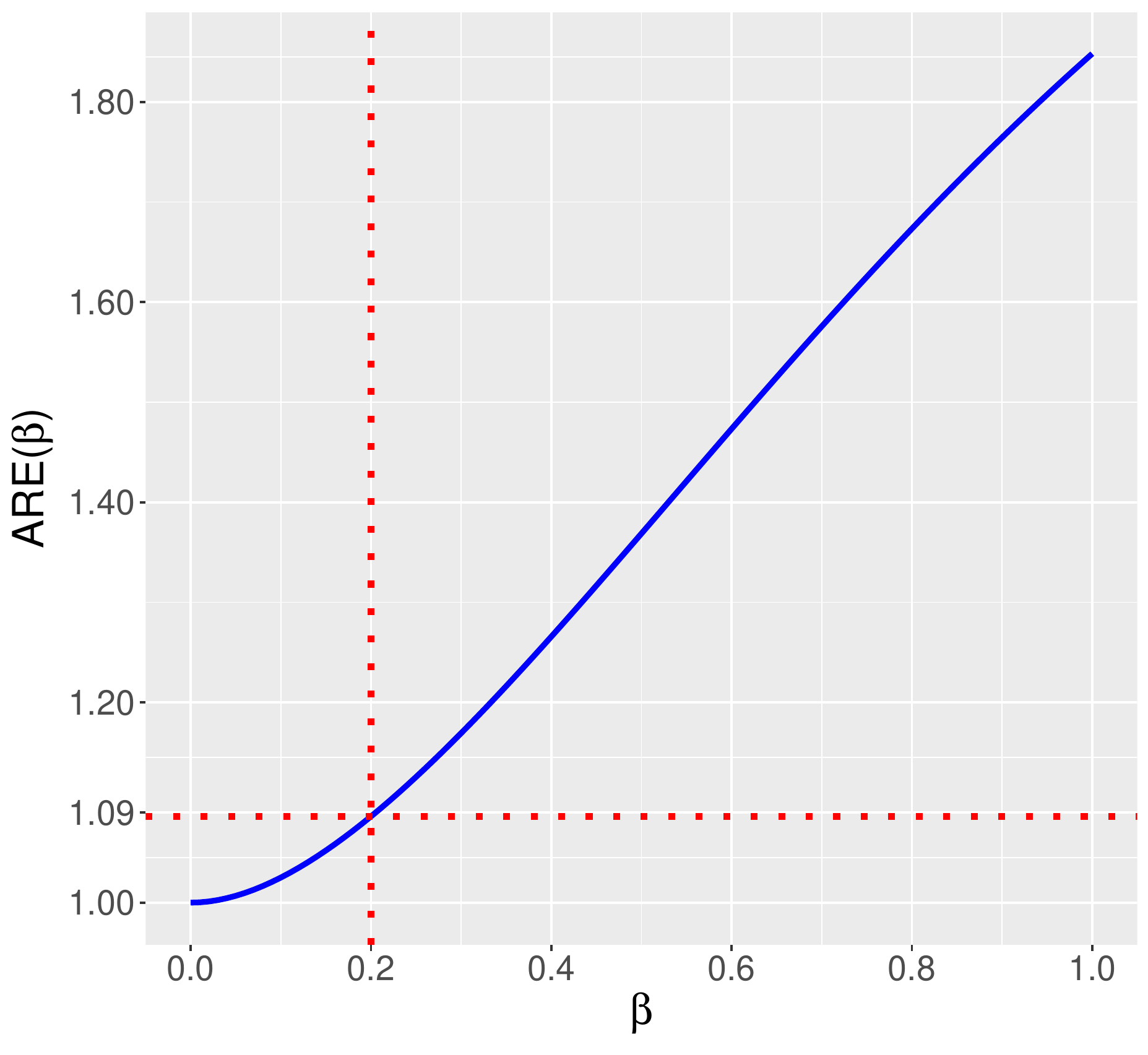}
	\caption{Asymptotic Relative
		Efficiency for the Breusch-Pagan $\beta$-score LM test-statistic.}%
	\label{fig:AREplot}%
\end{figure}

The expression given in (\ref{ARE}), for pure normally distributed data, depends only on
$\beta$ (there is no dependence on $\alpha$, $\boldsymbol{\theta}_{0}$ or
$\boldsymbol{\delta}$), being actually increasing on $\beta$ (see Figure
\ref{fig:AREplot}). This means that, as it happens usually with the minimum DPD estimators,
the most efficient test-statistics are the ones closer in the value of $\beta$
to the classical one, $\beta=0$. It is observed that the ARE is $1.09$ at
$\beta=0.2$, hence for $\beta=0.2$, $R_{n,\beta=0.2}%
(\widetilde{\boldsymbol{\theta}}_{\beta=0.2})$ needs about $9\%$ additional
observations to get the same power as $R_{n,\beta=0}%
(\widetilde{\boldsymbol{\theta}}_{\beta=0})$.

\begin{corollary}
\label{Corollary}The Breusch-Pagan $\beta $-score LM test-statistic for the
heteroscedastic linear model, (\ref{eq1}) or (\ref{eq2}), is given by%
\begin{equation}
	R_{n,\beta}(\widetilde{\boldsymbol{\theta}}_{\beta})=\frac{1}{\tfrac
		{2(2\beta^{2}+1)}{(2\beta+1)^{5/2}}-\frac{\beta^{2}}{(\beta+1)^{3}}%
	}\boldsymbol{v}^{T}(\widetilde{\boldsymbol{\theta}}_{\beta})\mathbb{\breve{Z}%
	}_{n}(\mathbb{\breve{Z}}_{n}^{T}\mathbb{\breve{Z}}_{n})^{-1}\mathbb{\breve{Z}%
	}_{n}^{T}\boldsymbol{v}(\widetilde{\boldsymbol{\theta}}_{\beta}%
	),\label{LM_TS3}%
\end{equation}
with $\boldsymbol{v}(\widetilde{\boldsymbol{\theta }}_{\beta })$ given by (%
\ref{vectorLM}). Once the two conditions established by (\ref{cond}) are verified,
the test-statistic, (\ref{LM_TS3}), is asymptotically a chi-square random
variable with $r$ degrees of freedom under homoscedastic null hypothesis.
\end{corollary}

\begin{proof}
	From Theorem \ref{ThM2b}, we can use (\ref{LM_TS2})\ rather than
	(\ref{LM_TS}). The matrix of the quadratic form to be inverted needs a
	previous block matrix inversion and selection of its $r\times r$ upper-left
	block, i.e.
	\[
	\left(  \boldsymbol{\bar{K}}_{n,\beta}^{-1}(\widetilde{\boldsymbol{\theta}%
	}_{\beta})\right)  _{11}=\frac{1}{\frac{2\beta^{2}+1}{2(2\beta+1)^{5/2}}%
		-\frac{\beta^{2}}{4(\beta+1)^{3}}}\boldsymbol{W}_{n,\beta,11}^{-1}%
	(\widetilde{\boldsymbol{\theta}}_{\beta}),
	\]
	with
	\begin{align*}
		\boldsymbol{W}_{n,\beta,11}^{-1}(\widetilde{\boldsymbol{\theta}}_{\beta}) &
		=\tfrac{n}{\widetilde{\sigma}^{4}(h^{\prime}(0))^{2}}\left(  \mathbb{Z}%
		_{n}^{T}\boldsymbol{H}_{n}\mathbb{Z}_{n}\right)  ^{-1},\\
		\boldsymbol{H}_{n} &  =\boldsymbol{I}_{n}-\tfrac{1}{n}\boldsymbol{1}%
		_{n}\boldsymbol{1}_{n}^{T},
	\end{align*}
	according to Lu and Shiou (2002, Theorem 2.1.-(ii)). From (\ref{vectorLM}) and
	(\ref{eq3b}), we have
	\begin{align*}
		R_{n,\beta}(\widetilde{\boldsymbol{\theta}}_{\beta}) &  =\tfrac{1}{4\left(
			\frac{2\beta^{2}+1}{2(2\beta+1)^{5/2}}-\frac{\beta^{2}}{4(\beta+1)^{3}%
			}\right)  }\boldsymbol{v}^{T}(\widetilde{\boldsymbol{\theta}}_{\beta
		})\mathbb{Z}_{n}\left(  \mathbb{Z}_{n}^{T}\boldsymbol{H}_{n}\mathbb{Z}%
		_{n}\right)  ^{-1}\mathbb{Z}_{n}^{T}\boldsymbol{v}%
		(\widetilde{\boldsymbol{\theta}}_{\beta})\\
		&  =\tfrac{1}{4\left(  \frac{2\beta^{2}+1}{2(2\beta+1)^{5/2}}-\frac{\beta^{2}%
			}{4(\beta+1)^{3}}\right)  }\boldsymbol{v}^{T}(\widetilde{\boldsymbol{\theta}%
		}_{\beta})\left(  \boldsymbol{H}_{n}\mathbb{Z}_{n}\left(  \mathbb{Z}_{n}%
		^{T}\boldsymbol{H}_{n}\mathbb{Z}_{n}\right)  ^{-1}\mathbb{Z}_{n}%
		^{T}\boldsymbol{H}_{n}+\tfrac{1}{n}\boldsymbol{1}_{n}\boldsymbol{1}_{n}%
		^{T}\right)  \boldsymbol{v}(\widetilde{\boldsymbol{\theta}}_{\beta}).
	\end{align*}
	and the final expression (\ref{LM_TS3}) is followed from Theorem A.45 in\  Rao
	et al. (2008, page 503).
\end{proof}

\begin{remark}
	\label{interpretation}The classical Breusch-Pagan LM test-statistic ($\beta
	=0$),
	\begin{equation}
		R_{n}(\widetilde{\boldsymbol{\theta}})=\frac{1}{2}(\boldsymbol{g}%
		(\widetilde{\boldsymbol{\theta}})-\boldsymbol{1}_{r})^{T}\mathbb{\breve{Z}%
		}_{n}(\mathbb{\breve{Z}}_{n}^{T}\mathbb{\breve{Z}}_{n})^{-1}\mathbb{\breve{Z}%
		}_{n}^{T}(\boldsymbol{g}(\widetilde{\boldsymbol{\theta}})-\boldsymbol{1}%
		_{r}),\label{ClBP}%
	\end{equation}
	matches (5.17) of Godfrey (1989, page 128). This expression justifies the
	two-fold ordinary least squares (OLS) regression procedure for the computation
	of the classical Breusch-Pagan LM test-statistic. The most general
	interpretation for the new proposal is also easily deducted. The value of the
	Breusch-Pagan $\beta$-score LM test-statistic is based on a projection of
	$\boldsymbol{v}(\widetilde{\boldsymbol{\theta}}_{\beta})$, a vector
	constructed from the vector of squared standardized $\beta$-residuals,
	$\boldsymbol{g}(\widetilde{\boldsymbol{\theta}}_{\beta})$, taking
	$\boldsymbol{Y}$ as response to be adjusted in the linear regression model,
	with $\mathbb{X}_{n}$ as design matrix. Such a projection is on the vector
	space of dimension $r+1$, generated by the columns of $\mathbb{\breve{Z}}_{n}%
	$, whenever $(r+1)+(p+1)<n$. The practical procedure is as follows:
	\textquotedblleft The value of (\ref{LM_TS3}) can be calculated as the
	explained sum of squares from the OLS regression of the response
	$\boldsymbol{v}(\widetilde{\boldsymbol{\theta}}_{\beta})$ over the matrix of
	explanatory variables $\mathbb{\breve{Z}}_{n}$, multiplied by $(\tfrac
	{2(2\beta^{2}+1)}{(2\beta+1)^{5/2}}-\tfrac{\beta^{2}}{(\beta+1)^{3}})^{-1}%
	$\textquotedblright.
\end{remark}
\section{Influence Function Analysis\label{Sec4}}

The influence function approach was introduced in Hampel (1968, 1974) for
estimators, being in essence based on the Gateaux derivative of a
functional. It was fully developed in Hampel et al. (1986). It is the most
important tool to analyze the robustness of statistical procedures, due to
its adaptability to analyze either robustness of estimators or
test-statistics.

In this section we are going to introduce the second order influence
function of the $\beta $-score LM test-statistics for non-identically
distributed observation. The following scheme will be followed. We will
calculate first the influence function of the minimum DPD estimator of the
LM, under the null hypothesis, and based on it, later, the influence
function of the $\beta $-score LM test-statistics is obtained. To fully
justify the robustness of the test-statistic, in principle it is not enough
with analyzing the raw influence function, it requires to prove in addition
the stability of the significance level and power of the $\beta $-score LM
test-statistics\ under data contamination.

Let $\boldsymbol{T}_{n,\beta}(\underline{G})$ denote funcional associated with
$\widetilde{\boldsymbol{\theta}}_{2}$, where $\underline{G}=(G_{1}%
,\ldots,G_{n})^{T}$ is the vector of true distribution
observation-by-observation in the sample, $Y_{i}$, $i=1,\ldots,n$. Under the
assumption that the true distribution belongs to the parametric model
associated with the homoscedastic linear regression, the same vector is
denoted by $\underline{F}_{\boldsymbol{\theta}_{0}}=(F_{1,\boldsymbol{\theta
	}_{0}},\ldots,F_{n,\boldsymbol{\theta}_{0}})^{T}$.

\begin{theorem}
\label{Th_IF}The $\beta $-score LM tests for non-identically distributed
individuals, has the usual first-order influence function equals zero, and
its second-order influence function, or self standardized IF of $\boldsymbol{%
\lambda }_{n,\beta }(\widetilde{\boldsymbol{\theta }}_{\beta })$, is given by%
\begin{equation}
\mathcal{IF}^{(2)}\left( \boldsymbol{y},R_{n,\beta }\left( \boldsymbol{%
\theta }_{0}\right) ,\underline{F}_{\boldsymbol{\theta }_{0}}\right) =\frac{1%
}{n}\sum_{i=1}^{n}\mathcal{IF}_{i}^{(2)}\left( y_{i},R_{n,\beta }\left( 
\boldsymbol{\theta }_{0}\right) ,\underline{F}_{\boldsymbol{\theta }%
_{0}}\right) ,  \label{IF2}
\end{equation}%
where 
\begin{align}
& \mathcal{IF}_{i}^{(2)}\left( y_{i},R_{n,\beta }\left( \boldsymbol{\theta }%
_{0}\right) ,\underline{F}_{\boldsymbol{\theta }_{0}}\right)  \label{IF2b} \\
& =\boldsymbol{u}_{i,\beta }^{T}\left( y_{i};\boldsymbol{\theta }\right) 
\boldsymbol{J}_{\beta }^{-1}(\boldsymbol{\theta })\boldsymbol{M}(\boldsymbol{%
\theta }_{0})\left[ \boldsymbol{M}^{T}(\boldsymbol{\theta }_{0})\boldsymbol{V%
}_{\beta }(\boldsymbol{\theta }_{0})\boldsymbol{M}(\boldsymbol{\theta }_{0})%
\right] ^{-1}\boldsymbol{M}^{T}(\boldsymbol{\theta })\boldsymbol{J}_{\beta
}^{-1}(\boldsymbol{\theta })\boldsymbol{u}_{i,\beta }\left( y_{i};%
\boldsymbol{\theta }\right) .  \notag
\end{align}
\end{theorem}

Taking into account that $\boldsymbol{u}_{i,\beta }\left( y_{i};\boldsymbol{%
\theta }_{0}\right) $ is bounded, for any $\beta >0$, then it is fulfilled
that for any $\beta >0$, $\mathcal{IF}^{(2)}\left( \boldsymbol{y},R_{n,\beta
}\left( \boldsymbol{\theta }_{0}\right) ,\underline{F}_{\boldsymbol{\theta }%
_{0}}\right) $ is bounded with respect to the $i$-th observation for $%
i=1,\ldots ,n$, and the higher the value of $\beta >0$, the outliers suffer
from a greater down-weighting effect. For the classic score tests ($\beta =0$%
), $\mathcal{IF}^{(2)}\left( \boldsymbol{y},R_{n,\beta =0}\left( \boldsymbol{%
\theta }_{0}\right) ,\underline{F}_{\boldsymbol{\theta }_{0}}\right) $ is
unbounded because $\boldsymbol{u}_{i,\beta =0}\left( y_{i};\boldsymbol{%
\theta }_{0}\right) =\boldsymbol{s}_{i,\boldsymbol{\theta }_{0}}(y_{i})$ is
unbounded. The influence function of all the observations is bounded if only
if the influence functions associated with pairs of individuals are all
bounded, i.e. the influence function of all the observations is bounded if
only if $\beta >0$. The gross-error sensitivity (GES) of $R_{n,\beta }\left( 
\boldsymbol{\theta }_{0}\right) $ for a sample of non-identically
distributed observations is defined as%
\begin{equation*}
\mathcal{GES}(R_{n,\beta }\left( \boldsymbol{\theta }_{0}\right)
)=\sup_{i\in \{1,\ldots ,n\}}\mathcal{GES}_{i}(R_{n,\beta }\left( 
\boldsymbol{\theta }_{0}\right) )=\sup_{i\in \{1,\ldots ,n\}}\sup_{y_{i}\in 
\mathbb{R}
}\mathcal{IF}_{i}^{(2)}\left( y_{i},R_{n,\beta }\left( \boldsymbol{\theta }%
_{0}\right) ,\underline{F}_{\boldsymbol{\theta }_{0}}\right) ,
\end{equation*}%
where the GES of $R_{n,\beta }\left( \boldsymbol{\theta }_{0}\right) $ for
the $i$-th individual, $\mathcal{GES}_{i}(R_{n,\beta }\left( \boldsymbol{%
\theta }_{0}\right) )$, matches the square of the so-called self
standardized GES of $\boldsymbol{\lambda }_{n,\beta }(\widetilde{\boldsymbol{%
\theta }}_{\beta })$, $\boldsymbol{m}(\widetilde{\boldsymbol{\theta }}%
_{\beta })$, or $\widetilde{\boldsymbol{\theta }}_{\beta }$. A finite (an
infinite) GES implies a bounded (an unbounded) second order influence
function of $R_{n,\beta }\left( \boldsymbol{\theta }_{0}\right) $. Different
parameterizations yield a unique value of $R_{n,\beta =0}\left( \boldsymbol{%
\theta }_{0}\right) $, self standardized IF and self standardized\ GES
(introduced for the first time in Krasker and Welsch, 1982). This is the
so-called invariance property of the Rao-type test statistics, not fulfilled
by other commonly used test-statistics such as the likelihood ratio and Wald
type tests.

\begin{corollary}
In the particular case of $\boldsymbol{\theta }=(\boldsymbol{\theta }%
_{1}^{T},\boldsymbol{\theta }_{2}^{T})^{T}$ and $\Theta _{0}=\{\boldsymbol{%
\theta }\in \Theta :\boldsymbol{\theta }_{1}=\boldsymbol{\theta }_{1,0}\}$,
the $\beta $-score LM test for non-identically distributed individuals,
given in (\ref{LM_TS}), has the second-order influence function (\ref{IF2}),
where%
\begin{equation}
\mathcal{IF}_{i}^{(2)}\left( y_{i},R_{n,\beta }\left( \boldsymbol{\theta }%
_{0}\right) ,\underline{F}_{\boldsymbol{\theta }_{0}}\right) =\boldsymbol{u}%
_{i,\beta ,1}^{T}(\boldsymbol{\theta }_{0})\left( \boldsymbol{K}_{n,\beta
}^{-1}(\boldsymbol{\theta }_{0})\right) _{11}\boldsymbol{u}_{i,\beta ,1}(%
\boldsymbol{\theta }_{0}),  \label{IF2c}
\end{equation}%
and%
\begin{align*}
\boldsymbol{u}_{i,\beta ,1}\left( y_{i};\boldsymbol{\theta }\right) & =(%
\boldsymbol{I}_{r},\boldsymbol{0}_{r\times (p-r)})\boldsymbol{u}_{i,\beta
}\left( y_{i};\boldsymbol{\theta }\right) , \\
\left( \boldsymbol{K}_{n,\beta }^{-1}(\boldsymbol{\theta })\right) _{11}& =(%
\boldsymbol{I}_{r},\boldsymbol{0}_{r\times (p-r)})\boldsymbol{K}_{n,\beta
}^{-1}(\boldsymbol{\theta })(\boldsymbol{I}_{r},\boldsymbol{0}_{r\times
(p-r)})^{T}.
\end{align*}
\end{corollary}

Based on $\underline{F}_{\boldsymbol{\theta }_{0,n}(\boldsymbol{\delta })}(%
\boldsymbol{y})$, with $\boldsymbol{\theta }_{0,n}(\boldsymbol{\delta })\in
\Theta _{1,n}$, we shall consider the $\epsilon $-contaminated distribution
function under a sequence of local Pitman-type alternatives, with respect to
the degenerated distribution function, $\Lambda _{\boldsymbol{y}}$, 
\begin{equation*}
\underline{F}_{\boldsymbol{n,\epsilon }}(\boldsymbol{y};\boldsymbol{\theta }%
_{0},\boldsymbol{\delta })=\left( 1-\frac{\epsilon }{\sqrt{n}}\right) 
\underline{F}_{\boldsymbol{\theta }_{0,n}(\boldsymbol{\delta })}(\boldsymbol{%
y})+\frac{\epsilon }{\sqrt{n}}\Lambda _{\boldsymbol{y}}.
\end{equation*}

\begin{theorem}
\label{contaminatedPower}For the $\epsilon $-contaminated distribution
function and under a sequence of local Pitman-type alternatives $H_{1,n}:%
\boldsymbol{\theta }\in \Theta _{1,n}$, with $\Theta _{1,n}$ given by (\ref%
{pitman}), the $\beta $-score LM test-statistics for non-identically
distributed individuals, $R_{n,\beta }(\widetilde{\boldsymbol{\theta }}%
_{\beta ,\epsilon })$, is asymptotically a chi-square random variable with $%
r $ degrees of freedom an noncentrality parameter%
\begin{equation}
\nu _{\beta }(\boldsymbol{\theta }_{0},\boldsymbol{\delta }_{\beta ,\epsilon
})=\nu _{\beta }(\boldsymbol{\theta }_{0},\boldsymbol{\delta })+2\epsilon 
\boldsymbol{\delta }^{T}\mathcal{IF}\left( \boldsymbol{y},\boldsymbol{%
\lambda }_{n,\beta }\left( \boldsymbol{\theta }_{0}\right) ,\underline{F}_{%
\boldsymbol{\theta }_{0}}\right) +\epsilon ^{2}\mathcal{IF}^{(2)}\left( 
\boldsymbol{y},R_{n,\beta }\left( \boldsymbol{\theta }_{0}\right) ,%
\underline{F}_{\boldsymbol{\theta }_{0}}\right) ,  \label{NCP2}
\end{equation}%
with $\nu _{\beta }(\boldsymbol{\theta }_{0},\boldsymbol{\delta })$ given by
(\ref{NCP}),%
\begin{align}
\mathcal{IF}\left( \boldsymbol{y},\boldsymbol{\lambda }_{n,\beta }\left( 
\boldsymbol{\theta }_{0}\right) ,\underline{F}_{\boldsymbol{\theta }%
_{0}}\right) & =\sum\limits_{i=1}^{n}\mathcal{IF}\left( y_{i},\boldsymbol{%
\lambda }_{n,\beta }\left( \boldsymbol{\theta }_{0}\right) ,\underline{F}_{%
\boldsymbol{\theta }_{0}}\right) ,  \label{IFLambda} \\
\mathcal{IF}\left( y_{i},\boldsymbol{\lambda }_{n,\beta }\left( \boldsymbol{%
\theta }_{0}\right) ,\underline{F}_{\boldsymbol{\theta }_{0}}\right) & =%
\left[ \boldsymbol{M}^{T}(\boldsymbol{\theta }_{0})\boldsymbol{V}_{\beta }(%
\boldsymbol{\theta }_{0})\boldsymbol{M}(\boldsymbol{\theta }_{0})\right]
^{-1}\mathcal{IF}\left( y_{i},\boldsymbol{m}\left( \boldsymbol{\theta }%
_{0}\right) ,\underline{F}_{\boldsymbol{\theta }_{0}}\right) ,  \notag \\
\mathcal{IF}\left( y_{i},\boldsymbol{m}\left( \boldsymbol{\theta }%
_{0}\right) ,\underline{F}_{\boldsymbol{\theta }_{0}}\right) & =\boldsymbol{M%
}^{T}(\boldsymbol{\theta }_{0})\mathcal{IF}\left( y_{i},\boldsymbol{\theta }%
_{0},\underline{F}_{\boldsymbol{\theta }_{0}}\right) ,  \notag \\
\mathcal{IF}\left( y_{i},\boldsymbol{\theta }_{0},\underline{F}_{\boldsymbol{%
\theta }_{0}}\right) & =\boldsymbol{J}_{\beta }^{-1}(\boldsymbol{\theta }%
_{0})\boldsymbol{u}_{i,\beta }\left( y_{i};\boldsymbol{\theta }_{0}\right) ,
\notag
\end{align}%
$\mathcal{IF}^{(2)}\left( \boldsymbol{y},R_{n,\beta }\left( \boldsymbol{%
\theta }_{0}\right) ,\underline{F}_{\boldsymbol{\theta }_{0}}\right) $ by (%
\ref{IF2}). Moreover, the corresponding asymptotic power function, given the
nominal level $\alpha $, is 
\begin{equation}
\pi _{\boldsymbol{\epsilon ,}\beta }(\boldsymbol{y};\boldsymbol{\theta }_{0},%
\boldsymbol{\delta })=\lim_{n\rightarrow \infty }P(R_{n,\beta }(\widetilde{%
\boldsymbol{\theta }}_{\beta })>\chi _{r,\alpha }^{2}|\underline{F}_{%
\boldsymbol{n,\epsilon }}(\boldsymbol{y};\boldsymbol{\theta }_{0},%
\boldsymbol{\delta }))=Q_{\frac{r}{2}}(\nu _{\beta }(\boldsymbol{\theta }%
_{0},\boldsymbol{\delta }_{\beta ,\epsilon }),\chi _{r,\alpha }^{2}),
\label{power2}
\end{equation}%
where $Q_{M}(a,b)$ is the Marcum $Q$-function.
\end{theorem}

The (asymptotic) power IF is defined as%
\begin{equation*}
\mathcal{PIF}_{\beta }(\boldsymbol{y};\boldsymbol{\theta }_{0},\boldsymbol{%
\delta })=\left. \frac{d}{d\epsilon }\pi _{\boldsymbol{\epsilon ,}\beta }(%
\boldsymbol{y};\boldsymbol{\theta }_{0},\boldsymbol{\delta })\right\vert
_{\epsilon =0}.
\end{equation*}

\begin{theorem}
\label{ThPIF}For the $\epsilon $-contaminated distribution function and
under a sequence of local Pitman-type alternatives $H_{1,n}:\boldsymbol{%
\theta }\in \Theta _{1,n}$, with $\Theta _{1,n}$ given by (\ref{pitman}),
the power inflluence function of the $\beta $-score LM test-statistics for
non-identically distributed individuals, $R_{n,\beta }(\widetilde{%
\boldsymbol{\theta }}_{\beta ,\epsilon })$,\ with nominal level $\alpha $,
is given by 
\begin{align*}
\mathcal{PIF}_{\beta }(\boldsymbol{y};\boldsymbol{\theta }_{0},\boldsymbol{%
\delta })& =2\nu _{\beta }(\boldsymbol{\theta }_{0},\boldsymbol{\delta })%
\boldsymbol{\delta }^{T}\mathcal{IF}\left( \boldsymbol{y},\boldsymbol{%
\lambda }_{n,\beta }\left( \boldsymbol{\theta }_{0}\right) ,\underline{F}_{%
\boldsymbol{\theta }_{0}}\right) \\
& \times \left[ Q_{\frac{r}{2}+1}(\nu _{\beta }(\boldsymbol{\theta }_{0},%
\boldsymbol{\delta }),\chi _{r,\alpha }^{2})-Q_{\frac{r}{2}}(\nu _{\beta }(%
\boldsymbol{\theta }_{0},\boldsymbol{\delta }),\chi _{r,\alpha }^{2})\right]
,
\end{align*}%
with $\mathcal{IF}\left( \boldsymbol{y},\boldsymbol{\lambda }_{n,\beta
}\left( \boldsymbol{\theta }_{0}\right) ,\underline{F}_{\boldsymbol{\theta }%
_{0}}\right) $ given by (\ref{IFLambda}).
\end{theorem}

We may observe that the boundedness of the PIF\ is equivalent to the
boundedness of the IF of the minimum DPDs of parameter $\boldsymbol{\theta }$%
. In fact, it is well-known, that the IF is bounded for $\beta >0$ and
unbounded for $\beta =0$.

Taking $\boldsymbol{\delta =0}_{r}$ and $\boldsymbol{\theta }=\boldsymbol{%
\theta }_{0}$, since $\boldsymbol{\theta }_{n}(\boldsymbol{0}_{r})=%
\boldsymbol{\theta }_{0}\in \Theta _{0}$, based on $\underline{F}_{%
\boldsymbol{\theta }_{0}}(\boldsymbol{y})$, we shall consider $\underline{F}%
_{\boldsymbol{n,\epsilon }}(\boldsymbol{y};\boldsymbol{0}_{r})$ to be the $%
\epsilon $-contaminated distribution function under the null hypothesis and 
\begin{equation*}
\alpha _{\boldsymbol{\epsilon ,}\beta }(\boldsymbol{y})=\pi _{\boldsymbol{%
\epsilon ,}\beta }(\boldsymbol{y};\boldsymbol{\theta }_{0},\boldsymbol{0}%
_{r})=\lim_{n\rightarrow \infty }P(R_{n,\beta }(\widetilde{\boldsymbol{%
\theta }}_{\beta })>\chi _{r,\alpha }^{2}|\underline{F}_{\boldsymbol{%
n,\epsilon }}(\boldsymbol{y};\boldsymbol{\theta }_{0})).
\end{equation*}%
In this setting, $\mathcal{LIF}_{\beta }(\boldsymbol{y})=\left. \frac{d}{%
d\epsilon }\alpha _{\boldsymbol{\epsilon ,}\beta }(\boldsymbol{y}%
)\right\vert _{\epsilon =0}=0$, which indicates a null IF of the
significance level.

\begin{example}
For the Breusch-Pagan $\beta $-score LM test-statistics, the second order
influence function is%
\begin{equation*}
\mathcal{IF}^{(2)}\left( \boldsymbol{y},R_{n,\beta }\left( \boldsymbol{0}%
_{r},\sigma _{0}^{2},\boldsymbol{\beta }_{0}\right) ,\underline{F}_{%
\boldsymbol{0}_{r},\sigma _{0}^{2},\boldsymbol{\beta }_{0}}\right) =\frac{1}{%
n}\sum_{i=1}^{n}\mathcal{IF}_{i}^{(2)}\left( y_{i},R_{n,\beta }\left( 
\boldsymbol{0}_{r},\sigma _{0}^{2},\boldsymbol{\beta }_{0}\right) ,%
\underline{F}_{\boldsymbol{0}_{r},\sigma _{0}^{2},\boldsymbol{\beta }%
_{0}}\right) ,
\end{equation*}%
where 
\begin{equation*}
\mathcal{IF}_{i}^{(2)}\left( y_{i},R_{n,\beta }\left( \boldsymbol{0}%
_{r},\sigma _{0}^{2},\boldsymbol{\beta }_{0}\right) ,\underline{F}_{%
\boldsymbol{0}_{r},\sigma _{0}^{2},\boldsymbol{\beta }_{0}}\right) =\frac{%
\exp \left\{ -\beta \frac{(y_{i}-\boldsymbol{x}_{i}^{T}\boldsymbol{\beta }%
_{0})^{2}}{\sigma _{0}^{2}}\right\} \left( \frac{(y_{i}-\boldsymbol{x}%
_{i}^{T}\boldsymbol{\beta }_{0})^{2}}{\sigma _{0}^{2}}-1\right) ^{2}}{%
4\left( \tfrac{2\beta ^{2}+1}{2(2\beta +1)^{5/2}}-\tfrac{\beta ^{2}}{4(\beta
+1)^{3}}\right) }\boldsymbol{z}_{i}^{T}\left( \mathbb{Z}_{n}^{T}\mathbb{Z}%
_{n}\right) ^{-1}\boldsymbol{z}_{i}.
\end{equation*}%
It is fulfilled that for any $\beta >0$, $\mathcal{IF}_{i}^{(2)}(y,R_{n,\beta
}(\sigma _{0}^{2},\boldsymbol{\beta }_{0}))$ is bounded with respect to $y$,
but unbounded for the classical Breusch-Pagan score LM test ($\beta =0$). In
addition, leverages do not affect negatively on the second order influence
function for the Breusch-Pagan $\beta $-score LM tests since the first
factor is bounded with respect to $\boldsymbol{x}_{i}$\ and the second
factor is also bounded with respect to $\boldsymbol{z}_{i}$. In particular,
for the simple linear regression, with $p=1$, $x_{i1}=z_{i1}=i$ ((\ref{cond}%
) is verified), $\beta _{0}=0$, $\beta _{1}=1$, $\alpha _{1}=0$, being $%
\underline{G}=\underline{F}_{\left( 0,\sigma ^{2},0,1\right) }$ the true
model, we get%
\begin{equation*}
\mathcal{IF}_{i}^{(2)}(y_{i},R_{n,\beta }\left( 0,\sigma ^{2},0,1\right) ,%
\underline{F}_{\left( 0,\sigma ^{2},0,1\right) })=\frac{\exp \left\{ -\beta 
\frac{(y_{i}-i)^{2}}{\sigma ^{2}}\right\} \left( \frac{(y_{i}-i)^{2}}{\sigma
^{2}}-1\right) ^{2}}{2\left( \tfrac{2\beta ^{2}+1}{2(2\beta +1)^{5/2}}-%
\tfrac{\beta ^{2}}{4(\beta +1)^{3}}\right) }\frac{3i^{2}}{(2n+1)(n^{2}-1)},
\end{equation*}%
In Figure \ref{fig} the second order IFs are plotted when $i=1$, $\sigma
^{2}\in \{2,6\}$,and $\beta \in \{0,0.15,0.2,0.3\}$. On the left and right
hand tails of the curves associated with $\beta >0$, the down-weighting
effect is clearly visualized while for $\beta =0$ the curve increases
indefinitely. The maximum of the curves is reached at a double solution on $%
y_{i}$ of $(y_{i}-i)^{2}=\frac{\beta +2}{\beta }\sigma ^{2}$, concluding 
\begin{equation}
\mathcal{GES}(R_{n,\beta }\left( 0,\sigma ^{2},0,1\right) )=6\left[ \beta
^{2}\exp \left\{ \beta +2\right\} \left( \tfrac{2\beta ^{2}+1}{2(2\beta
+1)^{5/2}}-\tfrac{\beta ^{2}}{4(\beta +1)^{3}}\right) \right] ^{-1}\frac{%
n^{2}}{(2n+1)(n^{2}-1)},  \label{GES}
\end{equation}%
a decreasing function on $\beta $, being finite only for $\beta >0$. Since $%
\mathcal{GES}(R_{n,\beta }\left( 0,\sigma ^{2},0,1\right) )=\infty $ for $%
\beta =0$, the stability of the classical Breusch-Pagan LM test-statistic
breaks down completely with outliers, while robustness of the Breusch-Pagan $%
\beta $-score LM test-statistic increases (GES decreases) as $\beta >0$
increases. 
\begin{figure}[tbp]
\par
\begin{center}
\subfloat[Case: $\sigma^2=2$]{{\includegraphics[trim=0cm 0cm 0cm 0cm,
clip=true,width=12cm]{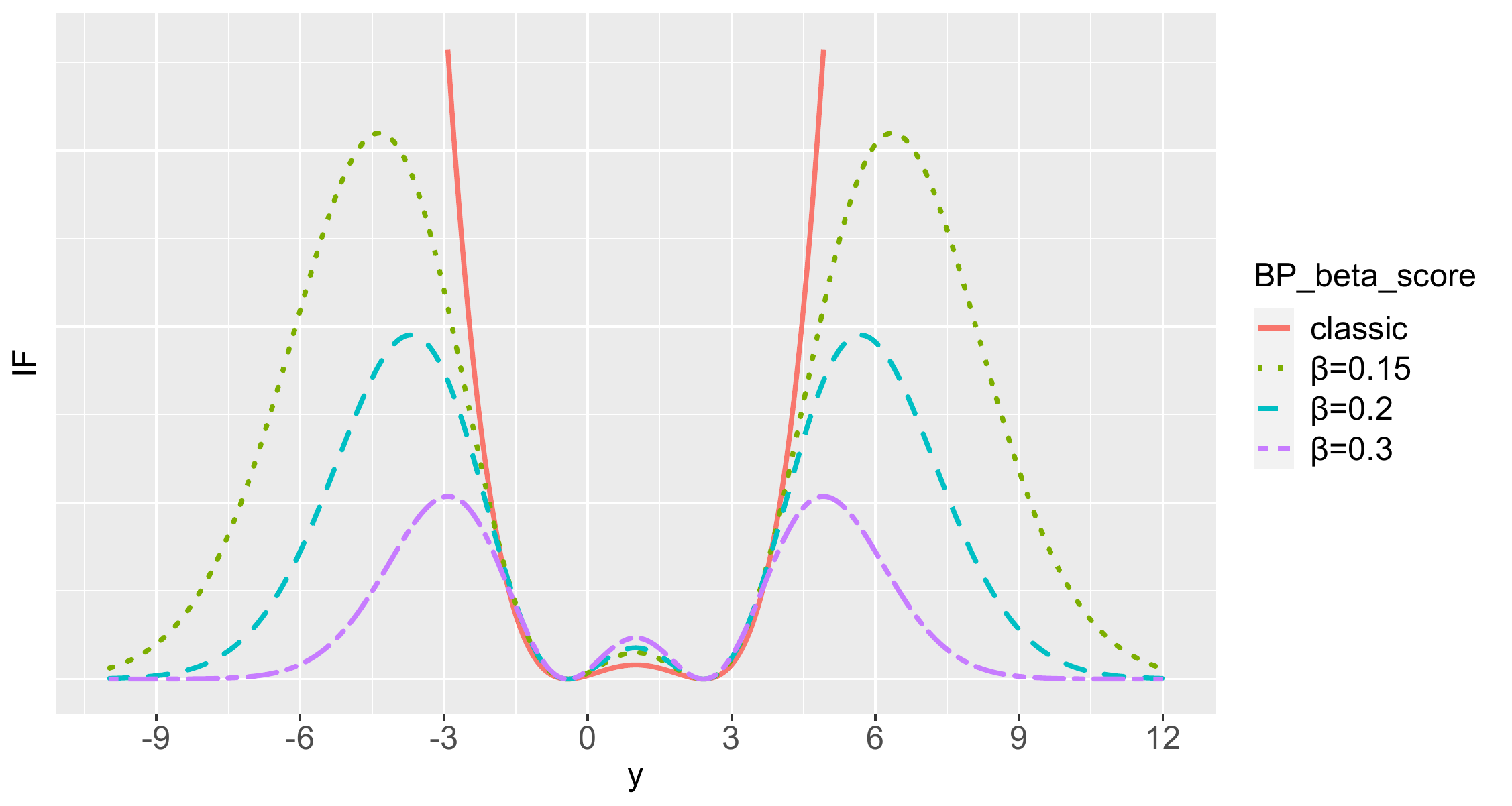} }} \qquad 
\subfloat[Case: $\sigma^2=6$]{{\includegraphics[trim=0cm 0cm 0cm 0cm,
clip=true,width=12cm]{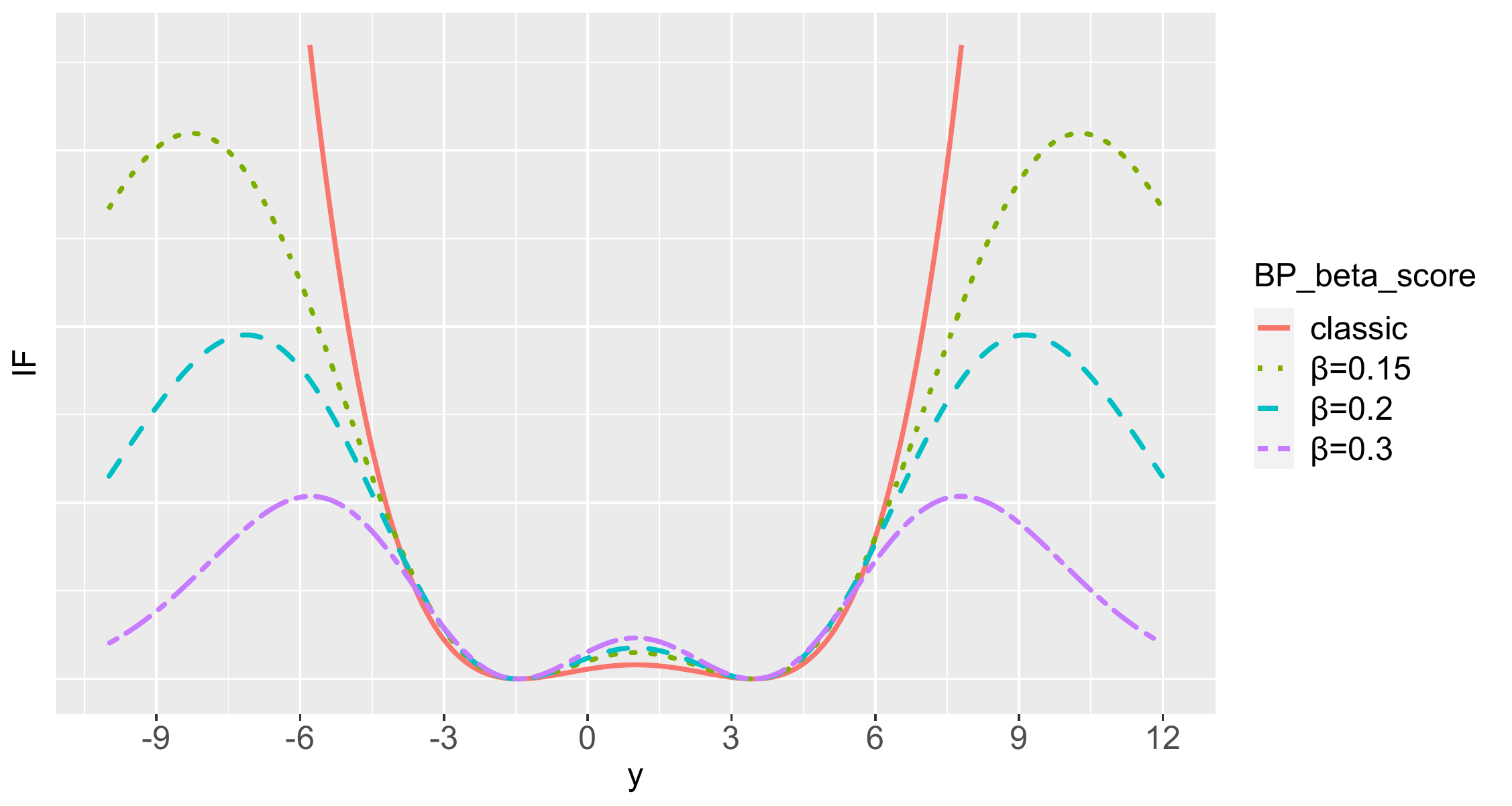} }}
\end{center}
\caption{$\mathcal{IF}_{i=1}^{(2)}(y_{i=1},R_{n,\protect\beta }\left( 0,%
\protect\sigma ^{2},0,1\right) ,\protect\underline{F}_{\left( 0,\protect%
\sigma ^{2},0,1\right) })$ for the Breusch-Pagan $\protect\beta $-score LM
tests of a simple linear regression.}
\label{fig}
\end{figure}
\end{example}

\section{Koenker's $\beta$-score test-statistic}

\label{Sec4b} Koenker (1981) extended the Breusch-Pagan test-statistic to make
it applicable also for non-normally distributed errors. Multiplying the
original Breusch-Pagan test-statistic by a constant, the new proposal was
found to have the same asymptotic distribution, with a new OLS interpretation
and relationship with the sample kurtosis. Similarly, we shall extend the
Breusch-Pagan $\beta$-score LM test in what we call the Koenker's $\beta
$-score test. For the classical Breusch-Pagan test ($\beta=0$), the average of
the cross-product of the estimating function is%
\[
\tfrac{1}{n}\boldsymbol{v}^{T}(\widetilde{\boldsymbol{\theta}}_{\beta
	=0})\boldsymbol{v}(\widetilde{\boldsymbol{\theta}}_{\beta=0})=\tfrac{1}%
{n}(\boldsymbol{g}(\widetilde{\boldsymbol{\theta}}_{\beta=0})-\boldsymbol{1}%
_{r})^{T}(\boldsymbol{g}(\widetilde{\boldsymbol{\theta}}_{\beta=0}%
)-\boldsymbol{1}_{r})=\tfrac{1}{n}\boldsymbol{g}^{T}%
(\widetilde{\boldsymbol{\theta}}_{\beta=0})\boldsymbol{g}%
(\widetilde{\boldsymbol{\theta}}_{\beta=0})-1,
\]
where $\tfrac{1}{n}\boldsymbol{g}^{T}(\widetilde{\boldsymbol{\theta}}%
_{\beta=0})\boldsymbol{g}(\widetilde{\boldsymbol{\theta}}_{\beta=0})$ is a
consistent estimator of the theoretical Kurtosis coefficient of the $i$-th
error (under the null hypothesis of homoscedasticity), i.e. $E\left[
g_{i}^{2}(\boldsymbol{\theta})\right]  =\frac{E[\epsilon_{i}^{4}%
	(\boldsymbol{\beta})]}{\sigma^{4}}$. In the genuine linear homoscedastic
regression, i.e. when $\epsilon_{i}(\boldsymbol{\beta})\sim\mathcal{N}%
(0,\sigma^{2})$, it holds the theoretical Kurtosis coefficient to be $3$.
Hence, it is almost straightforward to see that if we replace $2$ by
$\tfrac{1}{n}\boldsymbol{v}(\widetilde{\boldsymbol{\theta}}_{\beta
	=0})\boldsymbol{v}(\widetilde{\boldsymbol{\theta}}_{\beta=0})$ in the
denominator of (\ref{ClBP}), we get the Koenker's test-statistic%
\begin{equation}
	Q_{n}(\widetilde{\boldsymbol{\theta}}_{\beta=0})=\frac{1}{\tfrac{1}%
		{n}\boldsymbol{g}^{T}(\widetilde{\boldsymbol{\theta}}_{\beta=0})\boldsymbol{g}%
		(\widetilde{\boldsymbol{\theta}}_{\beta=0})-1}(\boldsymbol{g}%
	(\widetilde{\boldsymbol{\theta}}_{\beta=0})-\boldsymbol{1}_{r})^{T}%
	\mathbb{\breve{Z}}_{n}(\mathbb{\breve{Z}}_{n}^{T}\mathbb{\breve{Z}}_{n}%
	)^{-1}\mathbb{\breve{Z}}_{n}^{T}(\boldsymbol{g}(\widetilde{\boldsymbol{\theta
	}}_{\beta=0})-\boldsymbol{1}_{r}).\label{koenker_Cl}%
\end{equation}
In the same vein, the Koenker's $\beta$-score test-statistic is defined as%
\begin{equation}
	Q_{n}(\widetilde{\boldsymbol{\theta}}_{\beta})=\frac{1}{\tfrac{1}%
		{n}\boldsymbol{v}^{T}(\widetilde{\boldsymbol{\theta}}_{\beta})\boldsymbol{v}%
		(\widetilde{\boldsymbol{\theta}}_{\beta})}\boldsymbol{v}^{T}%
	(\widetilde{\boldsymbol{\theta}}_{\beta})\mathbb{\breve{Z}}_{n}(\mathbb{\breve
		{Z}}_{n}^{T}\mathbb{\breve{Z}}_{n})^{-1}\mathbb{\breve{Z}}_{n}^{T}%
	\boldsymbol{v}(\widetilde{\boldsymbol{\theta}}_{\beta}),\label{Koenker_beta}%
\end{equation}
where%
\[
\tfrac{1}{n}\boldsymbol{v}^{T}(\widetilde{\boldsymbol{\theta}}_{\beta
})\boldsymbol{v}(\widetilde{\boldsymbol{\theta}}_{\beta})=\tfrac{1}{n}\left(
\boldsymbol{g}(\widetilde{\boldsymbol{\theta}}_{\beta})-\boldsymbol{1}%
_{r}\right)  ^{T}\mathrm{diag}\left(  \exp(-\beta\boldsymbol{g}%
(\widetilde{\boldsymbol{\theta}}_{\beta}))\right)  \left(  \boldsymbol{g}%
(\widetilde{\boldsymbol{\theta}}_{\beta})-\boldsymbol{1}_{r}\right)
-\frac{\beta^{2}}{(\beta+1)^{3}}.
\]

\begin{corollary}
	\label{corKoenker}The Koenker's $\beta$-score test-statistic for the
	heteroscedastic linear model (\ref{eq1}) or (\ref{eq2}), given by
	(\ref{Koenker_beta}),  follows asymptotically a chi-square distribution with
	$r$ degrees of freedom.
\end{corollary}

\begin{proof}
	From the Weak Law of the Large Numbers (WLLN) and taking into account that
	$\widetilde{\boldsymbol{\theta}}_{\beta}\underset{n\rightarrow\infty
	}{\overset{p}{\longrightarrow}}\boldsymbol{\theta}$, it holds
	\begin{equation}
		\tfrac{1}{n}\left(  \boldsymbol{g}(\widetilde{\boldsymbol{\theta}}_{\beta
		})-\boldsymbol{1}_{r}\right)  ^{T}\mathrm{diag}\left(  \exp(-\beta
		\boldsymbol{g}(\widetilde{\boldsymbol{\theta}}_{\beta}))\right)  \left(
		\boldsymbol{g}(\widetilde{\boldsymbol{\theta}}_{\beta})-\boldsymbol{1}%
		_{r}\right)  \underset{n\rightarrow\infty}{\overset{p}{\longrightarrow}%
		}E\left[  \exp\{-\beta g_{i}(\boldsymbol{\theta})\}(g_{i}(\boldsymbol{\theta
		})-1)^{2}\right]  ,\label{cuadratic}%
	\end{equation}
	where $E\left[  \exp\{-\beta g_{i}(\boldsymbol{\theta})\}(g_{i}%
	(\boldsymbol{\theta})-1)^{2}\right]  =\varphi(-\beta)-2\frac{\partial
	}{\partial\beta}\varphi(-\beta)+\frac{\partial^{2}}{\partial\beta^{2}}%
	\varphi(-\beta)=\tfrac{2(2\beta^{2}+1)}{(2\beta+1)^{5/2}}$, with
	$\varphi(s)=(1-2s)^{-\frac{1}{2}}$ being the moment generating function of
	$g_{i}(\boldsymbol{\theta})\sim\chi_{1}^{2}$, valid among other values, for
	$s=-\beta<0$, and under normality of the errors. Taking into account that
	(\ref{cuadratic}) implies%
	\[
	\tfrac{1}{n}\boldsymbol{v}^{T}(\widetilde{\boldsymbol{\theta}}_{\beta
	})\boldsymbol{v}(\widetilde{\boldsymbol{\theta}}_{\beta}%
	)\underset{n\rightarrow\infty}{\overset{p}{\longrightarrow}}E\left[
	\exp\{-\beta g_{i}(\boldsymbol{\theta})\}(g_{i}(\boldsymbol{\theta}%
	)-1)^{2}\right]  -\frac{\beta^{2}}{(\beta+1)^{3}},
	\]
	with Corollary \ref{Corollary} and the Slutsky's theorem, the desired result
	is obtained.
\end{proof}

The calculation of (\ref{Koenker_beta}) can be made though a proper
interpretation of the least squared method as follows. It is $n$ times the
determination coefficient calculated from the linear regression of
$\boldsymbol{v}(\widetilde{\boldsymbol{\theta}}_{\beta})$, as response, over
explanatory variables given by $\mathbb{\breve{Z}}_{n}$. Hence, it remains
having the original appealing interpretation for the classical Koenker's
test-statistic. In Corollary \ref{corKoenker}, the validity of its asymptotic
distribution has been justified for normally distributed errors, in a similar
way done in Koenker (1981). Taking into account that the classical Koenker's
test, for non-normally distributed observations, was proven in a quantile
regression framework (see Koenker, 1982), the most general validity of the
asymptotic distribution of (\ref{Koenker_beta}) will be now justified. The
normality assumption of the error, (\ref{error}), of the linear regression model, (\ref{eq1}), can
be weakened to any centered distribution with finite variance, having second
order derivative of the moment generating function for the square of the
standardized errors.

\section{Example: Housing Price Data\label{Sec5}}

\begin{figure}[H]
	\centering
	\includegraphics[width=0.87\textwidth,trim=0 0 0 14mm,clip]{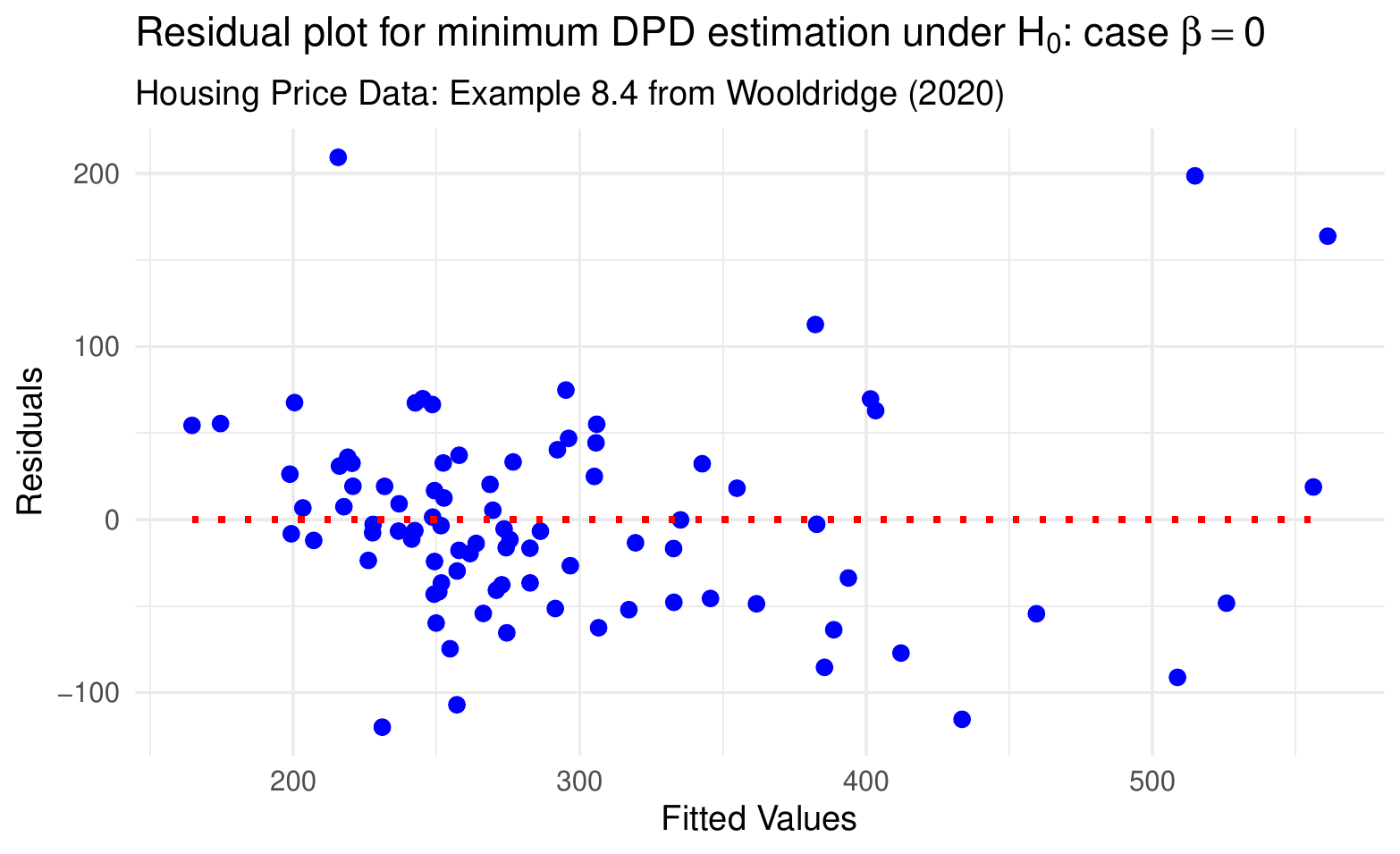}
	\caption{Residual plot for the Housing Price Data for estimation under OLS or MLEs ($\beta=0$).}
	\label{fig:Residual}
\end{figure}

Within the {\texttt{wooldridge}} {\texttt{R}}  package, the dataframe \texttt{hprice1} is related to a well-known example in Econometrics from Wooldridge (2020), which is referred to as the Housing Price Data. The dataset includes observations of 88 individuals and contains 10 variables, representing housing characteristics and sales information for houses sold in the Boston area, and it was collected from a publicly available database, published in 1990. We consider the heteroscedastic linear model (\ref{eq1}), \texttt{price }$=\beta _{0}+\beta _{1}$ \texttt{bdrms}$\ +\beta _{2}$ 
\texttt{lotsize}$\ +\beta _{3}$ \texttt{sqrft}, where \texttt{price} is the house price, in thousands of dollars ($y$), \texttt{bdrms}, the number of bedrooms ($x_1$), \texttt{lotsize}, the size of lot in square feet ($x_2$), \texttt{sqrft}, the size of house in square feet ($x_3$). Figure \ref{fig:Residual} suggests the presence of three outliers with the largest value of residuals, the observations with identification numbers $43$, $72$ and $76$ in the dataframe.

\begin{figure}[H]
	\centering
	\includegraphics[width=0.67\textwidth,clip]{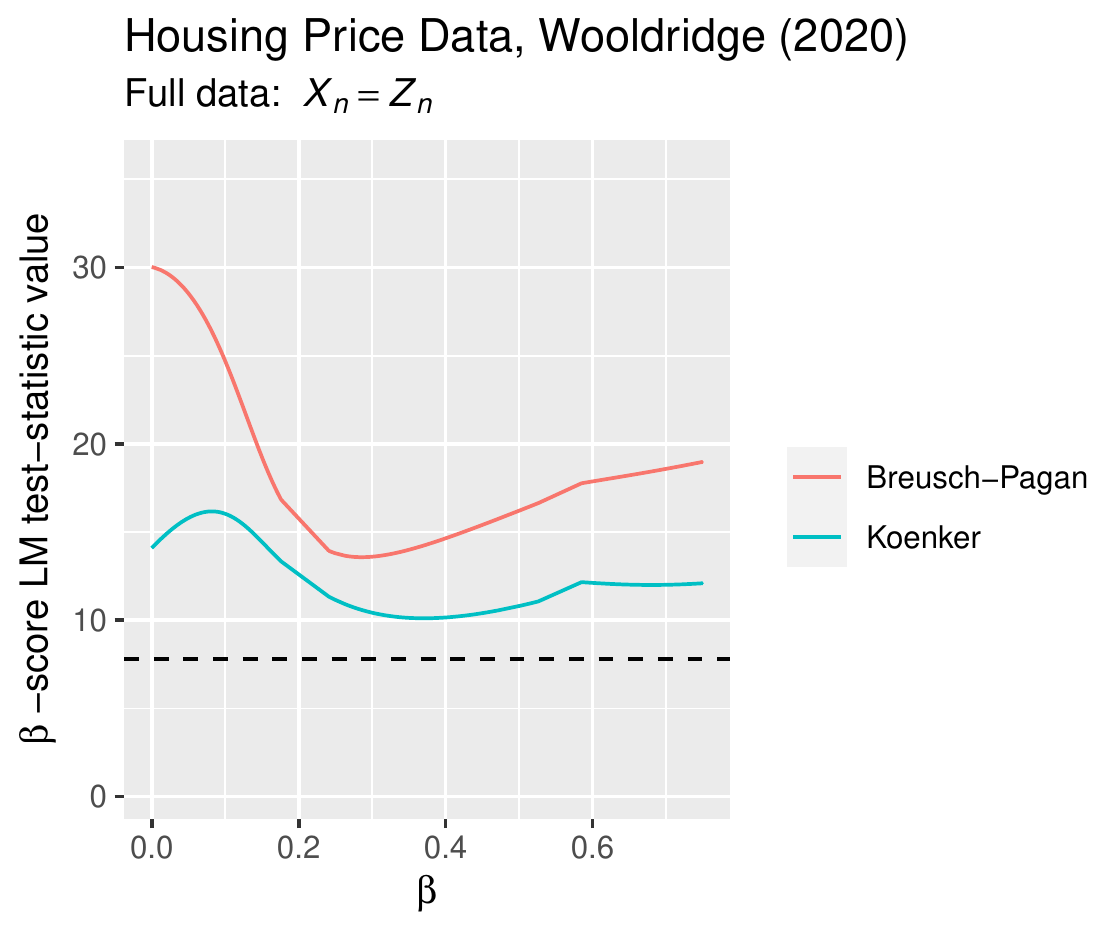}
	\hfill
	\includegraphics[width=0.67\textwidth,trim=0 0 0 7mm,clip]{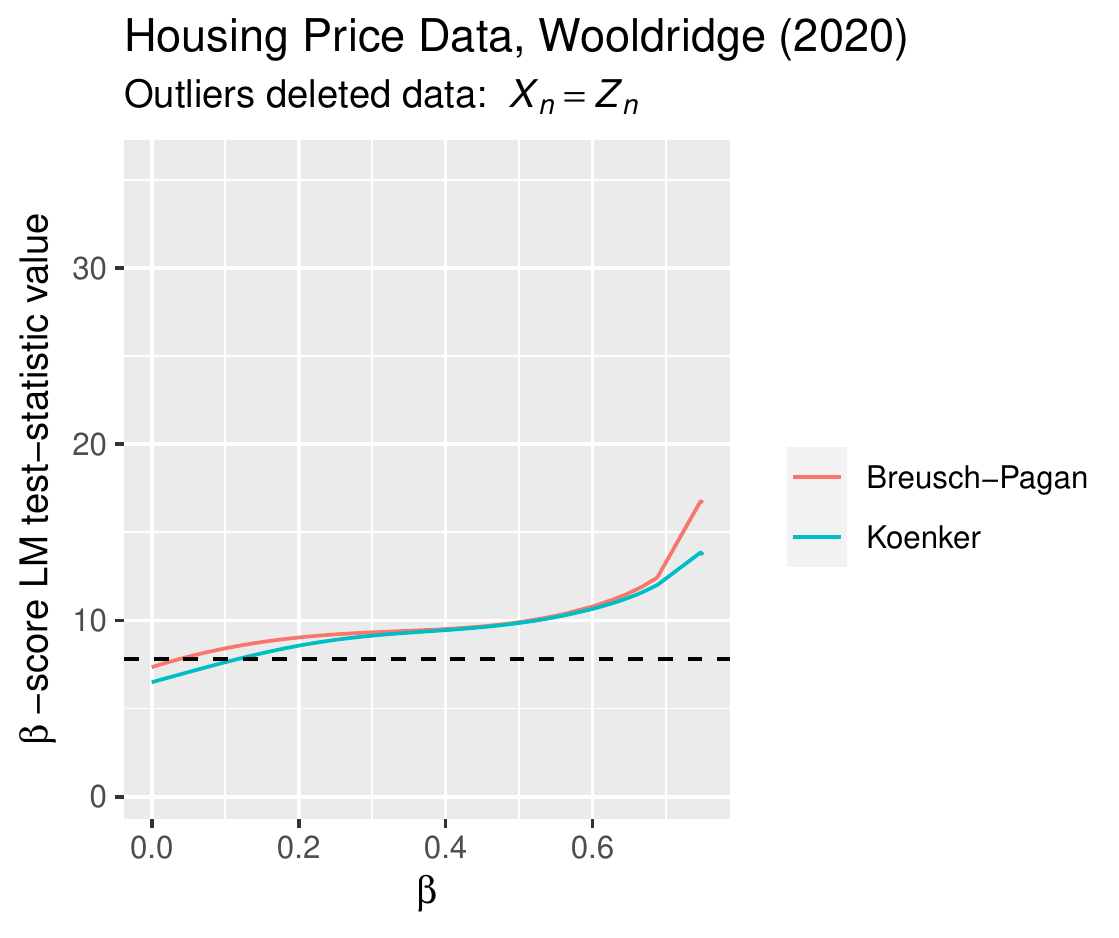}
\caption{Breusch-Pagan and Koenker's $\beta$-score LM test-statistics for the Housing Price Data and the case $\mathbb{X}_{n}=\mathbb{\breve{Z}}_{n}$.}
	\label{fig:WooldridgePlots1}
\end{figure}

\begin{figure}[H]
	\centering
\includegraphics[width=0.67\textwidth,clip]{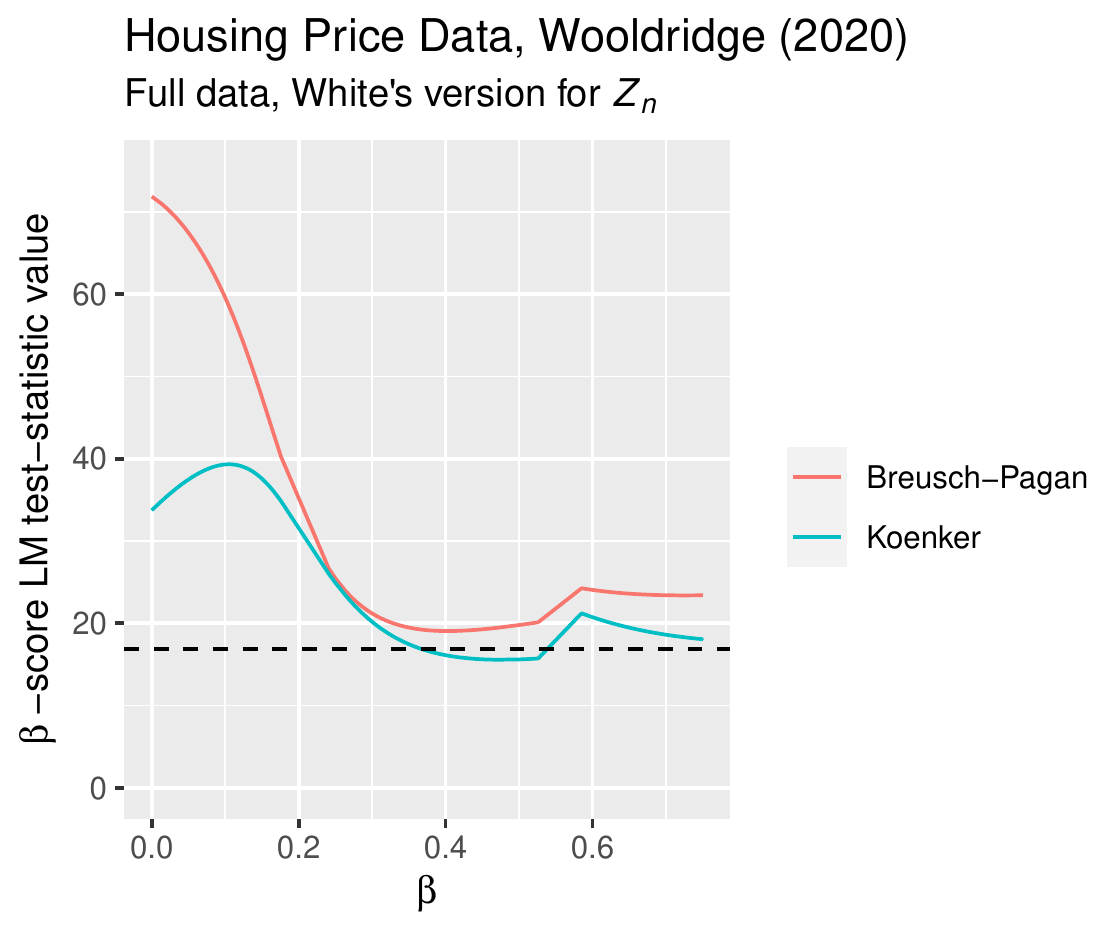}
\hfill
\includegraphics[width=0.67\textwidth,trim=0 0 0 7mm,clip]{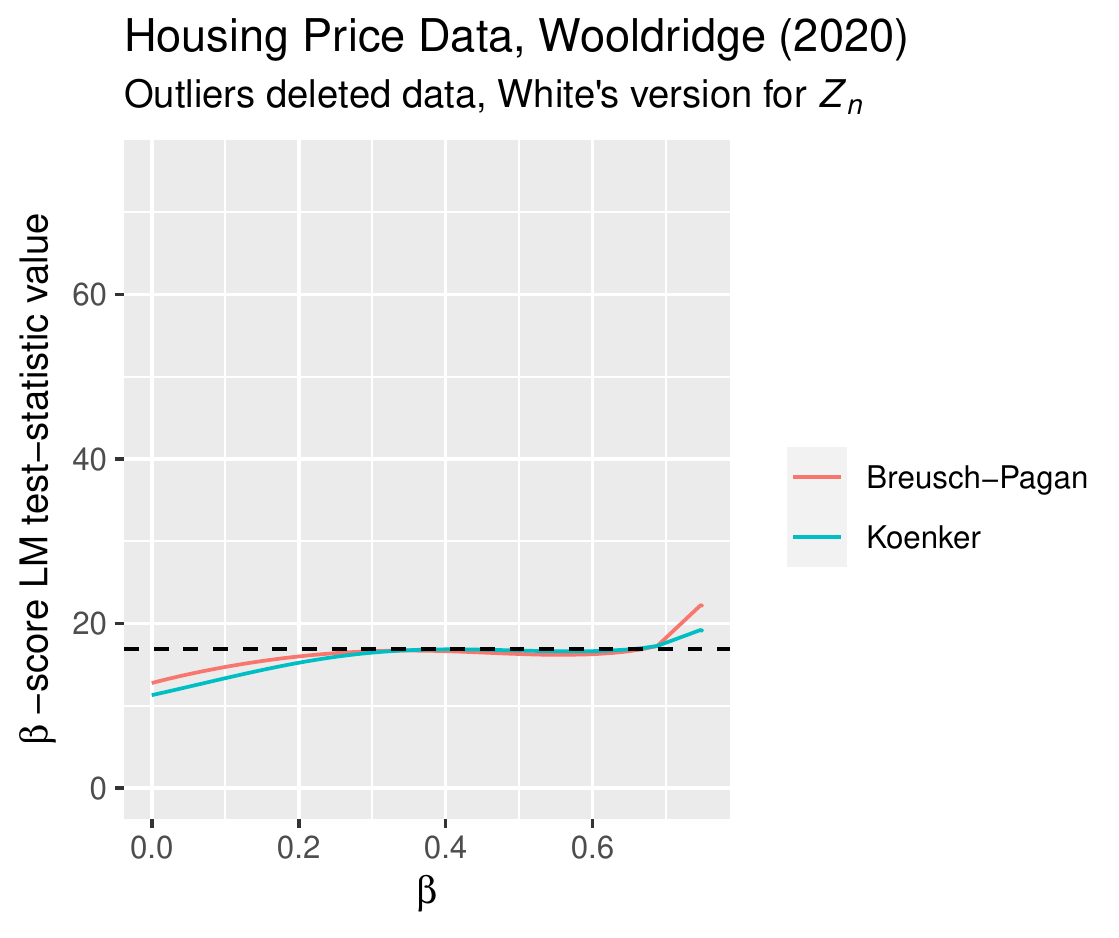}
	\caption{Breusch-Pagan and Koenker's $\beta$-score LM test-statistics for the Housing Price Data and White's version of $\mathbb{\breve{Z}}_{n}$.}
\label{fig:WooldridgePlots2}
\end{figure}

An extensive study of the performance of the Breusch-Pagan and Koenker's $\beta$-score LM test-statistics, $R_{n,\beta }(\widetilde{\boldsymbol{\theta }}_{\beta })$ and $Q_{n,\beta }(\widetilde{\boldsymbol{\theta }}_{\beta })$ respectively with tuning parameter interval $\beta \in [0, 0.75]$, has been performed taking in consideration either the full data or the outliers deleted data, dropping the three observations mentioned previously. We shall focus our attention on two different sources of heteroscedasticity depending on the choice of $\boldsymbol{\breve{z}}_{i}=(1,\boldsymbol{z}_{i})^{T}$, $i=1,\ldots,n$. The case  $\boldsymbol{x}_{i}=\boldsymbol{\breve{z}}_{i}$, $i=1,\ldots,n$ ($r=p=3$), is the one considered in Figure \ref{fig:WooldridgePlots1} and the left hand sides of Tables \ref{tab1:p.values}-\ref{tab2:p.values}, while in Figure \ref{fig:WooldridgePlots2} and the right hand sides of Tables \ref{tab1:p.values}-\ref{tab2:p.values} the White's additive heteroscedastic model is taken into account,%
\begin{equation*}
	h(\boldsymbol{z}_{i}^{T}\boldsymbol{\alpha })=1+\boldsymbol{z}_{i}^{T}%
	\boldsymbol{\alpha }=1+\sum_{j=1}^{3}x_{ij}\alpha
	_{j}+\sum_{j=1}^{3}\sum_{k\geq j}^{3}x_{ij}x_{ik}\alpha _{jk},
\end{equation*}%
with $\boldsymbol{z}_{i}=(%
\boldsymbol{x}_{i}^{T},\mathrm{vech}^{T}(\boldsymbol{x}_{i}\boldsymbol{x}%
_{i}^{T}))^{T}$, $vech(\boldsymbol{x}_{i}\boldsymbol{x}%
_{i})$ being the half vectorization of matrix $\boldsymbol{x}_{i}\boldsymbol{x}%
_{i}^{T}$, and $r= p+\frac{p(p+1)}{2}=9>p=3$. In the last section of the Appendix several computational details are given, useful for both estimation and test-statistics. Based on the $95\%$ quantile of the $\chi_{r}^{2}$ with $r$ being $3$ or $9$ depending on the choice of source of heteroscedasticity, the dashed line In Figures \ref{fig:WooldridgePlots1}-\ref{fig:WooldridgePlots2} is indicating the threshold for rejection of the test, considering the asymptotic distribution.

\afterpage{\clearpage}

The motivation behind this example is the suspicion of being made a non-appropriate decision with the classical tests, the ones with $\beta=0$. With both classical tests under consideration, the null hypothesis of homoscedasticity is clearly rejected, while it is accepted when deleting outliers. Our purpose is to provide further insight in the decision to be made, based on the study of the newly proposed test-statistics along the complete tuning parameter interval $[0,0.75]$. The strength of $R_{n,\beta }(\widetilde{\boldsymbol{\theta }}_{\beta })$ and $Q_{n,\beta }(\widetilde{\boldsymbol{\theta }}_{\beta })$, when $\beta>0$, is their damping effect in presence of outliers, and therefore avoiding the deletion of outliers. It is however, for illustrative purposes, interesting to analyze and compare the behaviour of the test-statistics when dropping outliers. 

In view of the plots shown in Figure \ref{fig:WooldridgePlots1} and the left hand sides of Tables \ref{tab1:p.values}-\ref{tab2:p.values}, we find enough evidence to reject the null hypothesis of homoscedasticity with $0.05$ significance level but in view of Figure \ref{fig:WooldridgePlots2} and the right hand sides of Tables \ref{tab1:p.values}-\ref{tab2:p.values} it is not so clear for both test-statistics since the Koenker's $\protect\beta $-score
LM test-statistic, when $\beta \in (0.36,0.53)$, is in the acceptance region. However, if we consider $0.056$ significance level, both test-statistics are in favour of the same decision of rejecting along the whole tuning parameter interval $[0.0.75]$.  This decision is not consistent with the one made in Berenguer-Rico and Wilms (2021). Indeed, the approach behind their proposal is different as well, since the method is constructed under the preference of dropping outliers from data.
\begin{table}[H]
	\centering 
	\begin{tabular}{ccccccc}
		\toprule & \multicolumn{3}{c}{case $\mathbb{X}_{n}=\mathbb{\breve{Z}}%
			_{n}$} & \multicolumn{3}{c}{White's version of $\mathbb{\breve{Z}}_{n}$} \\ 
		\midrule & $\beta =0$ & $\beta =0.3$ & $\beta =0.6$ & $\beta =0$ & $\beta
		=0.3$ &  $\beta
		=0.6$\\ 
		\cmidrule(lr){2-4}\cmidrule(l){5-7} $R_{n,\beta }(\widetilde{\boldsymbol{%
				\theta }}_{\beta })$ & 1.364e-06 & 3.501e-03 & 	4.6799e-04 & 6.559e-12 & 1.179e-02 & 0.0128
		\\ 
		$Q_{n,\beta }(\widetilde{\boldsymbol{\theta }}_{\beta })$ & 2.782e-03 & 
		1.526e-02 & 6.9919e-03 & 9.952e-05 & 1.684e-02 & 0.0137 \\ 
		\bottomrule &  &  &  &  &  & 
	\end{tabular}
	\caption{p-values of the Breusch-Pagan and Koenker's $\protect\beta $-score
		LM test-statistics, with full Housing Price Data.}
	\label{tab1:p.values}
\end{table}
\begin{table}[H]
	\centering 
	\begin{tabular}{ccccccc}
		\toprule & \multicolumn{3}{c}{case $\mathbb{X}_{n}=\mathbb{\breve{Z}}%
			_{n}$} & \multicolumn{3}{c}{White's version of $\mathbb{\breve{Z}}_{n}$} \\ 
		\midrule & $\beta =0$ & $\beta =0.3$ & $\beta =0.6$ & $\beta =0$ & $\beta
		=0.3$ & $\beta =0.6$ \\ 
		\cmidrule(lr){2-4}\cmidrule(l){5-7} $R_{n,\beta }(\widetilde{\boldsymbol{%
				\theta }}_{\beta })$ & 0.0615 & 0.0252 & 0.0041 & 0.1743 & 0.0550 &  0.0618\\ 
		$Q_{n,\beta }(\widetilde{\boldsymbol{\theta }}_{\beta })$ & 0.0898 & 0.0275
		& 0.0138 & 0.2576 & 0.0582 &  0.0548\\ 
		\bottomrule &  &  &  &  &  & 
	\end{tabular}
	\caption{p-values of the Breusch-Pagan and Koenker's $\protect\beta $-score
		LM test-statistics, with outliers deleted Housing Price Data.}
	\label{tab2:p.values}
\end{table}

\section{Concluding remarks\label{Sec6}}

This paper proposes Breusch-Pagan $\beta $-score LM tests for heterogeneity
in linear regression models. Prior to the proposal, the needed general
theory for its derivation has been developed over a framework of composite
null hypothesis with independent but non-identically distributed
observations, since Basu et al. (2022) did only consider identically
distributed observations. The Breusch-Pagan $\beta $-score LM
test-statistics depend on a tuning paramater $\beta \geq 0$ and the case $%
\beta =0$ is obtained as a limit of $\beta \rightarrow 0^{+}$ from the cases 
$\beta >0$. This family of test-statistics considers the classical
Breusch-Pagan LM test (Breusch and Pagan, 1979) as a particular case, being
just the case $\beta =0$. Efficiency and the robustness can be balanced
through the tuning parameter $\beta \geq 0$, being $\beta =0$ the most
efficient one, useful in case of having clean data free from outliers, while
robustness in the test-statistic is gained by increasing the value of $\beta
>0$, convenient in case of having outliers. A compromise between efficiency
and the robustness is reached by selecting properly an appropriate value of
the tuning parameter. 

Since Koenker (1981) published a modified version, robust against normality
assumption, this one remains being a reference of robustness, but our paper
fill the gap of proposing a specific test robust from outliers. The
challenge of generalizing our proposed $\beta $-score LM tests for
heterogeneity to be resistant from both simultaneously, normality assumption
and outliers, has been in addition addressed through the Koenker's $\beta $-score test.

When outliers are present in the data, in order to avoid results derived
from a biased and inefficient version of the classical Breusch-Godfried LM
test, based on Breusch (1978) and Godfrey (1978b), we are currently involved
in an ongoing paper which proposes $\beta $-score LM tests for testing
autocorrelation in the errors of a linear regression models. Apart from the aforementioned open problem, and based on the ideas
inherited from Breusch and Pagan (1980), the current paper could be the key
and reference for further challenging research on new applications of the $%
\beta $-score LM test for non-identically distributed individuals, useful in
Econometrics. In this regards, recently Halunga et al. (2017) have proposed
a robust version, against heteroscedasticity, of the Breusch-Pagan test for the null hypothesis of
zero cross-section correlation in dynamic panel data. Even though the model is
different and also the robust technique compared with our proposal, they comment \textquotedblleft If
the null is not rejected by the test, it would be more confidently concluded
that the rejection is not due to the heteroscedasticity, and OLS estimation
would be preferred. If the null is rejected, then a suitable estimation
procedure should be pursued\textquotedblright , which suits quite perfectly
in the context of the model in our paper and $\beta $-score LM tests could
be an attractive alternative methodology.

\appendix
\appendixpage

\section{Proof of Theorem \protect\ref{ThM}\label{AppA}}

We shall follow a similar scheme to Davidson and MacKinnon (2021, page 275),
in Section 8.9\ devoted to the classical LM test, but adapted to independent
and non-identical individuals and DPD based estimating functions. We may
start by considering the Taylor expansion of $\boldsymbol{u}_{i,\beta
}(Y_{i};\widetilde{\boldsymbol{\theta }}_{\beta })$ at the point $%
\boldsymbol{\theta }_{0}$\ for $i=1,\ldots ,n$%
\begin{equation*}
\boldsymbol{u}_{i,\beta }(Y_{i};\widetilde{\boldsymbol{\theta }}_{\beta })=%
\boldsymbol{u}_{i,\beta }(Y_{i};\boldsymbol{\theta }_{0})+\frac{\partial }{%
\partial \boldsymbol{\theta }}\left. \boldsymbol{u}_{i,\boldsymbol{\beta }%
}^{T}(Y_{i};\boldsymbol{\theta })\right\vert _{\boldsymbol{\theta =\theta }%
_{0}}(\widetilde{\boldsymbol{\theta }}_{\beta }-\boldsymbol{\theta }%
_{0})+o(\left\Vert \widetilde{\boldsymbol{\theta }}_{\beta }-\boldsymbol{%
\theta }_{0}\right\Vert ^{2}\boldsymbol{1}_{r}),
\end{equation*}%
we get 
\begin{equation*}
\sqrt{n}\boldsymbol{u}_{i,\beta }(Y_{i};\widetilde{\boldsymbol{\theta }}%
_{\beta })=\sqrt{n}\boldsymbol{u}_{i,\beta }(Y_{i};\boldsymbol{\theta }_{0})+%
\sqrt{n}\frac{\partial }{\partial \boldsymbol{\theta }}\left. \boldsymbol{u}%
_{i,\beta }^{T}(Y_{i};\boldsymbol{\theta })\right\vert _{\boldsymbol{\theta
=\theta }_{0}}(\widetilde{\boldsymbol{\theta }}_{\beta }-\boldsymbol{\theta }%
_{0})+o(\sqrt{n}\left\Vert \widetilde{\boldsymbol{\theta }}_{\beta }-%
\boldsymbol{\theta }_{0}\right\Vert ^{2}\boldsymbol{1}_{r}),
\end{equation*}%
This yields%
\begin{equation*}
\sqrt{n}\boldsymbol{U}_{n,\beta }(\widetilde{\boldsymbol{\theta }}_{\beta })=%
\sqrt{n}\boldsymbol{U}_{n,\beta }(\boldsymbol{\theta }_{0})+\frac{\partial }{%
\partial \boldsymbol{\theta }}\boldsymbol{U}_{n,\beta }^{T}(\boldsymbol{%
\theta }_{0})\sqrt{n}(\widetilde{\boldsymbol{\theta }}_{\beta }-\boldsymbol{%
\theta }_{0})+o(\sqrt{n}\left\Vert \widetilde{\boldsymbol{\theta }}_{\beta }-%
\boldsymbol{\theta }_{0}\right\Vert ^{2}\boldsymbol{1}_{r}).
\end{equation*}%
Adjusting the left hand side expression according to (\ref{eqA}) and taking
into account the Markov's WLLN%
\begin{equation*}
\frac{\partial }{\partial \boldsymbol{\theta }}\boldsymbol{U}_{n,\beta }^{T}(%
\boldsymbol{\theta }_{0})\overset{p}{\underset{n\rightarrow \infty }{%
\longrightarrow }}-\boldsymbol{J}_{\beta }(\boldsymbol{\theta }_{0}),
\end{equation*}%
it holds%
\begin{equation}
\sqrt{n}\boldsymbol{M}(\widetilde{\boldsymbol{\theta }}_{\beta })\boldsymbol{%
\lambda }_{n,\beta }(\widetilde{\boldsymbol{\theta }}_{\beta })=\sqrt{n}%
\boldsymbol{U}_{n,\beta }(\boldsymbol{\theta }_{0})-\boldsymbol{J}_{\beta }(%
\boldsymbol{\theta }_{0})\sqrt{n}(\widetilde{\boldsymbol{\theta }}_{\beta }-%
\boldsymbol{\theta }_{0})+o(\sqrt{n}\left\Vert \widetilde{\boldsymbol{\theta 
}}_{\beta }-\boldsymbol{\theta }_{0}\right\Vert ^{2}\boldsymbol{1}_{r}).
\label{eqAA}
\end{equation}%
On the other hand, from the Taylor expansion of $\boldsymbol{m}(\widetilde{%
\boldsymbol{\theta }}_{\beta })$ at the point $\boldsymbol{\theta }_{0}$, 
\begin{equation*}
\boldsymbol{m}(\widetilde{\boldsymbol{\theta }}_{\beta })=\boldsymbol{m}(%
\boldsymbol{\theta }_{0})+\boldsymbol{M}^{T}(\boldsymbol{\theta }_{0})(%
\widetilde{\boldsymbol{\theta }}_{\beta }-\boldsymbol{\theta }%
_{0})+o(\left\Vert \widetilde{\boldsymbol{\theta }}_{\beta }-\boldsymbol{%
\theta }_{0}\right\Vert ^{2}\boldsymbol{1}_{r}),
\end{equation*}%
and taking into account that(\ref{eqB}), we get%
\begin{equation}
\boldsymbol{0}_{r}=\boldsymbol{M}^{T}(\boldsymbol{\theta }_{0})\sqrt{n}(%
\widetilde{\boldsymbol{\theta }}_{\beta }-\boldsymbol{\theta }_{0})+o(\sqrt{n%
}\left\Vert \widetilde{\boldsymbol{\theta }}_{\beta }-\boldsymbol{\theta }%
_{0}\right\Vert ^{2}\boldsymbol{1}_{r}).  \label{eqBB}
\end{equation}%
Equations (\ref{eqAA})-(\ref{eqBB}) yield%
\begin{equation*}
\begin{bmatrix}
\boldsymbol{J}_{\beta }(\boldsymbol{\theta }_{0}) & \boldsymbol{M}(%
\boldsymbol{\theta }_{0}) \\ 
\boldsymbol{M}^{T}(\boldsymbol{\theta }_{0}) & \boldsymbol{0}_{p\times p}%
\end{bmatrix}%
\sqrt{n}%
\begin{bmatrix}
\widetilde{\boldsymbol{\theta }}_{\beta }-\boldsymbol{\theta }_{0} \\ 
\boldsymbol{\lambda }_{n,\beta }(\widetilde{\boldsymbol{\theta }}_{\beta })%
\end{bmatrix}%
=%
\begin{bmatrix}
\sqrt{n}\boldsymbol{U}_{n,\beta }(\boldsymbol{\theta }_{0}) \\ 
\boldsymbol{0}_{r}%
\end{bmatrix}%
+o(\sqrt{n}\left\Vert \widetilde{\boldsymbol{\theta }}_{\beta }-\boldsymbol{%
\theta }_{0}\right\Vert ^{2}\boldsymbol{1}_{r}).
\end{equation*}%
Hence,%
\begin{align}
\sqrt{n}%
\begin{bmatrix}
\widetilde{\boldsymbol{\theta }}_{\beta }-\boldsymbol{\theta }_{0} \\ 
\boldsymbol{\lambda }_{n,\beta }(\widetilde{\boldsymbol{\theta }}_{\beta })%
\end{bmatrix}%
& =%
\begin{bmatrix}
\boldsymbol{J}_{\beta }(\boldsymbol{\theta }_{0}) & \boldsymbol{M}(%
\boldsymbol{\theta }_{0}) \\ 
\boldsymbol{M}^{T}(\boldsymbol{\theta }_{0}) & \boldsymbol{0}_{r\times r}%
\end{bmatrix}%
^{-1}%
\begin{bmatrix}
\sqrt{n}\boldsymbol{U}_{n,\beta }(\boldsymbol{\theta }_{0}) \\ 
\boldsymbol{0}_{r}%
\end{bmatrix}%
+o(\sqrt{n}\left\Vert \widetilde{\boldsymbol{\theta }}_{\beta }-\boldsymbol{%
\theta }_{0}\right\Vert ^{2}\boldsymbol{1}_{r})  \notag \\
& =%
\begin{bmatrix}
\boldsymbol{\Sigma }_{11,\beta }(\boldsymbol{\theta }_{0})\sqrt{n}%
\boldsymbol{U}_{n,\beta }(\boldsymbol{\theta }_{0}) \\ 
\boldsymbol{\Sigma }_{21,\beta }(\boldsymbol{\theta }_{0})\sqrt{n}%
\boldsymbol{U}_{n,\beta }(\boldsymbol{\theta }_{0})%
\end{bmatrix}%
+o(\sqrt{n}\left\Vert \widetilde{\boldsymbol{\theta }}_{\beta }-\boldsymbol{%
\theta }_{0}\right\Vert ^{2}\boldsymbol{1}_{r}),  \label{r1}
\end{align}%
with%
\begin{equation*}
\begin{bmatrix}
\boldsymbol{J}_{\beta }(\boldsymbol{\theta }_{0}) & \boldsymbol{M}(%
\boldsymbol{\theta }_{0}) \\ 
\boldsymbol{M}^{T}(\boldsymbol{\theta }_{0}) & \boldsymbol{0}_{r\times r}%
\end{bmatrix}%
^{-1}=%
\begin{bmatrix}
\boldsymbol{\Sigma }_{11,\beta }(\boldsymbol{\theta }_{0}) & \boldsymbol{%
\Sigma }_{21,\beta }^{T}(\boldsymbol{\theta }_{0}) \\ 
\boldsymbol{\Sigma }_{21,\beta }(\boldsymbol{\theta }_{0}) & \boldsymbol{%
\Sigma }_{22,\beta }(\boldsymbol{\theta }_{0})%
\end{bmatrix}%
,
\end{equation*}%
where%
\begin{align*}
\boldsymbol{\Sigma }_{11,\beta }(\boldsymbol{\theta }_{0})& =\boldsymbol{J}%
_{\beta }^{-1}(\boldsymbol{\theta }_{0})+\boldsymbol{J}_{\beta }^{-1}(%
\boldsymbol{\theta }_{0})\boldsymbol{M}(\boldsymbol{\theta }_{0})\boldsymbol{%
\Sigma }_{22,\beta }(\boldsymbol{\theta }_{0})\boldsymbol{M}^{T}(\boldsymbol{%
\theta }_{0})\boldsymbol{J}_{\beta }^{-1}(\boldsymbol{\theta }_{0}) \\
& =\boldsymbol{J}_{\beta }^{-1}(\boldsymbol{\theta }_{0})-\boldsymbol{J}%
_{\beta }^{-1}(\boldsymbol{\theta }_{0})\boldsymbol{M}(\boldsymbol{\theta }%
_{0})(\boldsymbol{M}^{T}(\boldsymbol{\theta }_{0})\boldsymbol{J}_{\beta
}^{-1}(\boldsymbol{\theta }_{0})\boldsymbol{M}(\boldsymbol{\theta }%
_{0}))^{-1}\boldsymbol{M}^{T}(\boldsymbol{\theta }_{0})\boldsymbol{J}_{\beta
}^{-1}(\boldsymbol{\theta }_{0}), \\
\boldsymbol{\Sigma }_{22,\beta }(\boldsymbol{\theta }_{0})& =-(\boldsymbol{M}%
^{T}(\boldsymbol{\theta }_{0})\boldsymbol{J}_{\beta }^{-1}(\boldsymbol{%
\theta }_{0})\boldsymbol{M}(\boldsymbol{\theta }_{0}))^{-1}, \\
\boldsymbol{\Sigma }_{21,\beta }(\boldsymbol{\theta }_{0})& =-\boldsymbol{%
\Sigma }_{22,\beta }(\boldsymbol{\theta }_{0})\boldsymbol{M}^{T}(\boldsymbol{%
\theta }_{0})\boldsymbol{J}_{\beta }^{-1}(\boldsymbol{\theta }_{0}) \\
& =(\boldsymbol{M}^{T}(\boldsymbol{\theta }_{0})\boldsymbol{J}_{\beta }^{-1}(%
\boldsymbol{\theta }_{0})\boldsymbol{M}(\boldsymbol{\theta }_{0}))^{-1}%
\boldsymbol{M}^{T}(\boldsymbol{\theta }_{0})\boldsymbol{J}_{\beta }^{-1}(%
\boldsymbol{\theta }_{0}),
\end{align*}%
and under the existence assumption of (\ref{Jmatrix})-(\ref{Kmatrix}) it is
also verified the the following the multivariate version of the Lindeberg CLT for independent non-identically
distributed observations,%
\begin{equation}
\sqrt{n}\boldsymbol{U}_{n,\beta }(\boldsymbol{\theta }_{0})\overset{\mathcal{%
L}}{\underset{n\rightarrow \infty }{\longrightarrow }}\mathcal{N}(%
\boldsymbol{0}_{p},\boldsymbol{K}_{\beta }(\boldsymbol{\theta })).
\label{r2}
\end{equation}%
In view of the expression given in Definition \ref{def}\ for $\boldsymbol{%
\bar{K}}_{n,\beta }(\boldsymbol{\theta })$, this is just the estimator of
the variance for a centered sample, which is\ justified as follows%
\begin{align*}
\boldsymbol{\bar{K}}_{n,\beta }(\boldsymbol{\theta })& =n\mathrm{E}\left[ 
\boldsymbol{U}_{n,\beta }(\boldsymbol{\theta })\boldsymbol{U}_{n,\beta }^{T}(%
\boldsymbol{\theta })\right] \\
& =\frac{1}{n}\sum_{i=1}^{n}\sum_{j=1,j\neq i}^{n}\mathrm{E}[\boldsymbol{u}%
_{i,\beta }(Y_{i};\boldsymbol{\theta })]\mathrm{E}[\boldsymbol{u}_{j,%
\boldsymbol{\theta }}^{T}(Y_{i};\boldsymbol{\theta })]+\frac{1}{n}%
\sum_{i=1}^{n}\mathrm{E}[\boldsymbol{u}_{i,\beta }(Y_{i};\boldsymbol{\theta }%
)\boldsymbol{u}_{i,\beta }^{T}(Y_{i};\boldsymbol{\theta })] \\
& =\frac{1}{n}\sum_{i=1}^{n}\boldsymbol{K}_{i,\beta }\left( \boldsymbol{%
\theta }\right) ,
\end{align*}%
where $\mathrm{E}[\boldsymbol{u}_{i,\beta }(Y_{i};\boldsymbol{\theta })]=%
\boldsymbol{0}_{p}$, $i=1,\ldots ,n$. The expression $\boldsymbol{\bar{K}}%
_{n,\beta }(\boldsymbol{\theta })=\boldsymbol{\bar{J}}_{n,2\beta }(%
\boldsymbol{\theta })-\frac{1}{n}\sum_{i=1}^{n}\boldsymbol{\xi }_{i,\beta
}\left( \boldsymbol{\theta }\right) \boldsymbol{\xi }_{i,\beta }^{T}\left( 
\boldsymbol{\theta }\right) $, is straightforward from 
\begin{equation*}
\boldsymbol{K}_{i,\beta }\left( \boldsymbol{\theta }\right) =\boldsymbol{J}%
_{i,2\beta }(\boldsymbol{\theta })-\boldsymbol{\xi }_{i,\beta }\left( 
\boldsymbol{\theta }\right) \boldsymbol{\xi }_{i,\beta }^{T}\left( 
\boldsymbol{\theta }\right) ,\quad i=1,\ldots ,n,
\end{equation*}%
established in Basu et al. (1998) for identically distributed observations.
Finally, from (\ref{r1})-(\ref{r2}) the desired result is obtained.
\section{Proof of Theorem \protect\ref{ThM2}\label{AppB}}
From (\ref{r2}) we obtain%
\begin{align*}
	&  \sqrt{n}\boldsymbol{\lambda}_{n,\beta}(\widetilde{\boldsymbol{\theta}%
	}_{\beta})\underset{n\rightarrow\infty}{\overset{\mathcal{L}}{\longrightarrow
	}}\mathcal{N}\left(  \boldsymbol{0}_{r},\boldsymbol{\Sigma}_{21,\beta
	}(\boldsymbol{\theta}_{0})\boldsymbol{K}(\boldsymbol{\theta}_{0}%
	)\boldsymbol{\Sigma}_{21,\beta}^{T}(\boldsymbol{\theta}_{0})\right)  ,\\
	&  \left(  \boldsymbol{\bar{\Sigma}}_{n,21,\beta}^{T}%
	(\widetilde{\boldsymbol{\theta}}_{\beta})\boldsymbol{\bar{K}}_{n,\beta
	}(\widetilde{\boldsymbol{\theta}}_{\beta})\boldsymbol{\bar{\Sigma}%
	}_{n,21,\beta}^{T}(\widetilde{\boldsymbol{\theta}}_{\beta})\right)
	^{-\frac{1}{2}}\sqrt{n}\boldsymbol{\lambda}_{n,\beta}%
	(\widetilde{\boldsymbol{\theta}}_{\beta})\underset{n\rightarrow\infty
	}{\overset{\mathcal{L}}{\longrightarrow}}\mathcal{N}\left(  \boldsymbol{0}%
	_{r},\boldsymbol{I}_{r}\right)  ,\\
	&  R_{n,\beta}(\widetilde{\boldsymbol{\theta}}_{\beta})=n\boldsymbol{\lambda
	}_{n,\beta}^{T}(\widetilde{\boldsymbol{\theta}}_{\beta})\left(
	\boldsymbol{\bar{\Sigma}}_{n,22,\beta}(\widetilde{\boldsymbol{\theta}}_{\beta
	})\boldsymbol{\bar{K}}_{n,\beta}(\widetilde{\boldsymbol{\theta}}_{\beta
	})\boldsymbol{\bar{\Sigma}}_{n,22,\beta}^{T}(\widetilde{\boldsymbol{\theta}%
	}_{\beta})\right)  ^{-1}\boldsymbol{\lambda}_{n,\beta}%
	(\widetilde{\boldsymbol{\theta}}_{\beta})\underset{n\rightarrow\infty
	}{\overset{\mathcal{L}}{\longrightarrow}}\chi_{r}^{2}.
\end{align*}
\section{Proof of Theorem \protect\ref{contiguous}\label{AppE0}}

Let $\boldsymbol{\theta }_{n}$, the solution of $\sqrt{n}\boldsymbol{m}(%
\boldsymbol{\theta }_{n})=\boldsymbol{\delta }$ for any $n\in 
\mathbb{N}
$. Following the same scheme of the proof for Theorem \ref{AppA} but
replacing the role of $\boldsymbol{\theta }_{0}$\ by the one of $\boldsymbol{%
	\theta }_{n}$, we get%
\begin{equation*}
	\boldsymbol{m}(\widetilde{\boldsymbol{\theta }}_{\beta })=\boldsymbol{m}(%
	\boldsymbol{\theta }_{n})+\boldsymbol{M}^{T}(\boldsymbol{\theta }_{0})(%
	\widetilde{\boldsymbol{\theta }}_{\beta }-\boldsymbol{\theta }%
	_{n})+o(\left\Vert \widetilde{\boldsymbol{\theta }}_{\beta }-\boldsymbol{%
		\theta }_{n}\right\Vert ^{2}\boldsymbol{1}_{r}),
\end{equation*}%
and hence%
\begin{equation}
	\sqrt{n}\boldsymbol{m}(\widetilde{\boldsymbol{\theta }}_{\beta })=\sqrt{n}%
	\boldsymbol{m}(\boldsymbol{\theta }_{n})+\boldsymbol{M}^{T}(\boldsymbol{%
		\theta }_{0})\sqrt{n}(\widetilde{\boldsymbol{\theta }}_{\beta }-\boldsymbol{%
		\theta }_{n})+o(\sqrt{n}\left\Vert \widetilde{\boldsymbol{\theta }}_{\beta }-%
	\boldsymbol{\theta }_{n}\right\Vert ^{2}\boldsymbol{1}_{r}),  \label{mn}
\end{equation}%
where%
\begin{equation*}
	\sqrt{n}\boldsymbol{m}(\widetilde{\boldsymbol{\theta }}_{\beta })=\sqrt{n}%
	\boldsymbol{m}(\boldsymbol{\theta }_{0})+o_{p}(\boldsymbol{1}_{r}),
\end{equation*}%
and then%
\begin{align}
	\sqrt{n}%
	\begin{bmatrix}
		\widetilde{\boldsymbol{\theta }}_{\beta }-\boldsymbol{\theta }_{n} \\ 
		\boldsymbol{\lambda }_{n,\beta }(\widetilde{\boldsymbol{\theta }}_{\beta })%
	\end{bmatrix}%
	& =%
	\begin{bmatrix}
		\boldsymbol{J}_{\beta }(\boldsymbol{\theta }_{0}) & \boldsymbol{M}(%
		\boldsymbol{\theta }_{0}) \\ 
		\boldsymbol{M}^{T}(\boldsymbol{\theta }_{0}) & \boldsymbol{0}_{r\times r}%
	\end{bmatrix}%
	^{-1}%
	\begin{bmatrix}
		\sqrt{n}\boldsymbol{U}_{n,\beta }(\boldsymbol{\theta }_{0}) \\ 
		-\boldsymbol{\delta }%
	\end{bmatrix}%
	+o(\sqrt{n}\left\Vert \widetilde{\boldsymbol{\theta }}_{\beta }-\boldsymbol{%
		\theta }_{n}\right\Vert ^{2}\boldsymbol{1}_{r})  \notag \\
	\sqrt{n}\boldsymbol{\lambda }_{n,\beta }(\widetilde{\boldsymbol{\theta }}%
	_{\beta })& =\boldsymbol{\Sigma }_{21,\beta }(\boldsymbol{\theta }_{0})\sqrt{%
		n}\boldsymbol{U}_{n,\beta }(\boldsymbol{\theta }_{0})-\boldsymbol{\Sigma }%
	_{22,\beta }(\boldsymbol{\theta }_{0})\boldsymbol{\delta }+o(\sqrt{n}%
	\left\Vert \widetilde{\boldsymbol{\theta }}_{\beta }-\boldsymbol{\theta }%
	_{n}\right\Vert ^{2}\boldsymbol{1}_{r}).  \label{desar}
\end{align}%
Finally, from (\ref{r2}) we obtain%
\begin{align*}
	& \sqrt{n}\boldsymbol{\lambda }_{n,\beta }(\widetilde{\boldsymbol{\theta }}%
	_{\beta })\underset{n\rightarrow \infty }{\overset{\mathcal{L}}{%
			\longrightarrow }}\mathcal{N}\left( -\boldsymbol{\Sigma }_{22,\beta }(%
	\boldsymbol{\theta }_{0})\boldsymbol{\delta },\boldsymbol{\Sigma }_{21,\beta
	}(\boldsymbol{\theta }_{0})\boldsymbol{K}(\boldsymbol{\theta }_{0})%
	\boldsymbol{\Sigma }_{21,\beta }^{T}(\boldsymbol{\theta }_{0})\right) , \\
	& \left( \boldsymbol{\bar{\Sigma}}_{n,21,\beta }^{T}(\widetilde{\boldsymbol{%
			\theta }}_{\beta })\boldsymbol{\bar{K}}_{n,\beta }(\widetilde{\boldsymbol{%
			\theta }}_{\beta })\boldsymbol{\bar{\Sigma}}_{n,21,\beta }^{T}(\widetilde{%
		\boldsymbol{\theta }}_{\beta })\right) ^{-\frac{1}{2}}\sqrt{n}\boldsymbol{%
		\lambda }_{n,\beta }(\widetilde{\boldsymbol{\theta }}_{\beta })\underset{%
		n\rightarrow \infty }{\overset{\mathcal{L}}{\longrightarrow }}\mathcal{N}%
	\left( \boldsymbol{\mu }_{\beta }(\boldsymbol{\delta }),\boldsymbol{I}%
	_{r}\right) , \\
	& R_{n,\beta }(\widetilde{\boldsymbol{\theta }}_{\beta })=n\boldsymbol{%
		\lambda }_{n,\beta }^{T}(\widetilde{\boldsymbol{\theta }}_{\beta })\left( 
	\boldsymbol{\bar{\Sigma}}_{n,22,\beta }(\widetilde{\boldsymbol{\theta }}%
	_{\beta })\boldsymbol{\bar{K}}_{n,\beta }(\widetilde{\boldsymbol{\theta }}%
	_{\beta })\boldsymbol{\bar{\Sigma}}_{n,22,\beta }^{T}(\widetilde{\boldsymbol{%
			\theta }}_{\beta })\right) ^{-1}\boldsymbol{\lambda }_{n,\beta }(\widetilde{%
		\boldsymbol{\theta }}_{\beta })\underset{n\rightarrow \infty }{\overset{%
			\mathcal{L}}{\longrightarrow }}\chi _{r}^{2}(\nu _{\beta }(\boldsymbol{%
		\theta }_{0},\boldsymbol{\delta })),
\end{align*}%
with%
\begin{align*}
	\boldsymbol{\mu }_{\beta }(\boldsymbol{\delta })& =-\left( \boldsymbol{%
		\Sigma }_{21,\beta }(\boldsymbol{\theta }_{0})\boldsymbol{K}(\boldsymbol{%
		\theta }_{0})\boldsymbol{\Sigma }_{21,\beta }^{T}(\boldsymbol{\theta }%
	_{0})\right) ^{-\frac{1}{2}}\boldsymbol{\Sigma }_{22,\beta }(\boldsymbol{%
		\theta }_{0})\boldsymbol{\delta }=\boldsymbol{E}_{\beta }^{\frac{1}{2}}(%
	\boldsymbol{\theta })\boldsymbol{\delta }, \\
	\nu _{\beta }(\boldsymbol{\theta },\boldsymbol{\delta })& =\boldsymbol{%
		\delta }^{T}\boldsymbol{E}_{\beta }\boldsymbol{\delta },
\end{align*}%
and 
\begin{equation*}
	\boldsymbol{E}_{\beta }(\boldsymbol{\theta })=\left( \boldsymbol{M}^{T}(%
	\boldsymbol{\theta })\boldsymbol{J}_{\beta }^{-1}(\boldsymbol{\theta })%
	\boldsymbol{K}(\boldsymbol{\theta })\boldsymbol{J}_{\beta }^{-1}(\boldsymbol{%
		\theta })\boldsymbol{M}(\boldsymbol{\theta })\right) ^{-1}
\end{equation*}%
is the efficacy of the test.

\section{Proof of Theorem \protect\ref{Th_IF}\label{AppD}}

From Theorem \ref{ThM} it is fulfilled%
\begin{equation*}
\sqrt{n}\boldsymbol{\lambda }_{n,\beta }(\widetilde{\boldsymbol{\theta }}%
_{\beta })\underset{n\rightarrow \infty }{\overset{\mathcal{L}}{%
\longrightarrow }}\mathcal{N}\left( \boldsymbol{0}_{p},\boldsymbol{\Sigma }%
_{21,\beta }(\boldsymbol{\theta }_{0})\boldsymbol{K}_{\beta }(\boldsymbol{%
\theta }_{0})\boldsymbol{\Sigma }_{21,\beta }^{T}(\boldsymbol{\theta }%
_{0})\right) ,
\end{equation*}%
with 
\begin{align*}
\boldsymbol{\Sigma }_{21,\beta }(\boldsymbol{\theta }_{0})\boldsymbol{K}%
_{\beta }(\boldsymbol{\theta }_{0})\boldsymbol{\Sigma }_{21,\beta }^{T}(%
\boldsymbol{\theta }_{0})& =\boldsymbol{\Sigma }_{21,\beta }(\boldsymbol{%
\theta }_{0})\mathrm{E}\left[ \boldsymbol{U}_{n,\beta }\left( \boldsymbol{%
\theta }_{0}\right) \boldsymbol{U}_{n,\beta }^{T}\left( \boldsymbol{\theta }%
_{0}\right) \right] \boldsymbol{\Sigma }_{21,\beta }^{T}(\boldsymbol{\theta }%
_{0}) \\
& =\mathrm{E}\left[ \mathcal{IF}\left( \boldsymbol{y},\boldsymbol{\lambda }%
_{n,\beta }\left( \boldsymbol{T}_{n,\beta }(\underline{F}_{\boldsymbol{%
\theta }_{2}})\right) \right) \mathcal{IF}_{i}^{T}\left( \boldsymbol{y},%
\boldsymbol{\lambda }_{n,\beta }\left( \boldsymbol{T}_{n,\beta }(\underline{F%
}_{\boldsymbol{\theta }_{2}})\right) \right) \right] ,
\end{align*}%
$\boldsymbol{y}=(y_{1},\ldots ,y_{n})^{T}$. Hence,

The influence function associated with all the individuals for the Lagrange
multiplier, is given by a function $%
\mathbb{R}
^{n}\longrightarrow 
\mathbb{R}
$ such that%
\begin{equation*}
\mathcal{IF}\left( \boldsymbol{y},\boldsymbol{\lambda }_{n,\beta }\left( 
\boldsymbol{T}_{n,\beta }(\underline{F}_{\boldsymbol{\theta }_{2}})\right)
\right) =\frac{1}{n}\sum_{i=1}^{n}\mathcal{IF}_{i}\left( y_{i},\boldsymbol{%
\lambda }_{n,\beta }\left( \boldsymbol{T}_{n,\beta }(\underline{F}_{%
\boldsymbol{\theta }_{2}})\right) \right) ,
\end{equation*}%
where%
\begin{equation}
\mathcal{IF}_{i}\left( y_{i},\boldsymbol{\lambda }_{n,\beta }\left( 
\boldsymbol{T}_{n,\beta }(\underline{F}_{\boldsymbol{\theta }_{2}})\right)
\right) =\boldsymbol{\Sigma }_{21,\beta }(\boldsymbol{\theta }_{0})%
\boldsymbol{u}_{i,\beta }\left( y_{i};\boldsymbol{T}_{n,\beta }(\underline{F}%
_{\boldsymbol{\theta }_{2}})\right) ,\quad i=1,\ldots ,n;  \label{IF1}
\end{equation}%
i.e.%
\begin{equation*}
\mathcal{IF}\left( \boldsymbol{y},\boldsymbol{\lambda }_{n,\beta }\left( 
\boldsymbol{T}_{n,\beta }(\underline{F}_{\boldsymbol{\theta }_{2}})\right)
\right) =\boldsymbol{\Sigma }_{21,\beta }(\boldsymbol{\theta }_{0})%
\boldsymbol{U}_{n,\beta }\left( \boldsymbol{T}_{n,\beta }(\underline{F}_{%
\boldsymbol{\theta }_{2}})\right) .
\end{equation*}

\begin{align*}
& \mathcal{IF}_{i}^{(2)}\left( y_{i},R_{n,\beta }\left( \boldsymbol{T}%
_{n,\beta }(\underline{F}_{\boldsymbol{\theta }_{2}})\right) \right) = \\
& =\mathcal{IF}_{i}^{T}\left( y_{i},\boldsymbol{\lambda }_{n,\beta }\left( 
\boldsymbol{T}_{n,\beta }(\underline{F}_{\boldsymbol{\theta }_{2}})\right)
\right) \left[ \boldsymbol{\Sigma }_{21,\beta }(\boldsymbol{\theta }_{0})%
\boldsymbol{K}_{\beta }(\boldsymbol{\theta }_{0})\boldsymbol{\Sigma }%
_{21,\beta }^{T}(\boldsymbol{\theta }_{0})\right] ^{-1}\mathcal{IF}%
_{i}\left( y_{i},\boldsymbol{\lambda }_{n,\beta }\left( \boldsymbol{T}%
_{n,\beta }(\underline{F}_{\boldsymbol{\theta }_{2}})\right) \right) \\
& =\boldsymbol{u}_{i,\beta }^{T}\left( y_{i};\boldsymbol{T}_{n,\beta }(%
\underline{F}_{\boldsymbol{\theta }_{2}})\right) \boldsymbol{J}_{\beta
}^{-1}(\boldsymbol{\theta }_{0})\boldsymbol{M}(\boldsymbol{\theta }_{0})%
\left[ \boldsymbol{M}^{T}(\boldsymbol{\theta }_{0})\boldsymbol{J}_{\beta
}^{-1}(\boldsymbol{\theta }_{0})\boldsymbol{K}_{\beta }(\boldsymbol{\theta }%
_{0})\boldsymbol{J}_{\beta }^{-1}(\boldsymbol{\theta }_{0})\boldsymbol{M}(%
\boldsymbol{\theta }_{0})\right] ^{-1} \\
& \times \boldsymbol{M}^{T}(\boldsymbol{\theta }_{0})\boldsymbol{J}_{\beta
}^{-1}(\boldsymbol{\theta }_{0})\boldsymbol{u}_{i,\beta }\left( y_{i};%
\boldsymbol{T}_{n,\beta }(\underline{F}_{\boldsymbol{\theta }_{2}})\right) \\
& =\left\Vert \boldsymbol{M}^{T}(\boldsymbol{\theta }_{0})\boldsymbol{J}%
_{\beta }^{-1}(\boldsymbol{\theta }_{0})\boldsymbol{u}_{i,\beta }\left(
y_{i};\boldsymbol{T}_{n,\beta }(\underline{F}_{\boldsymbol{\theta }%
_{2}})\right) \right\Vert _{\boldsymbol{M}^{T}(\boldsymbol{\theta }_{0})%
\boldsymbol{J}_{\beta }^{-1}(\boldsymbol{\theta }_{0})\boldsymbol{K}_{\beta
}(\boldsymbol{\theta }_{0})\boldsymbol{J}_{\beta }^{-1}(\boldsymbol{\theta }%
_{0})\boldsymbol{M}(\boldsymbol{\theta }_{0})}^{2}.
\end{align*}%
\begin{align*}
& \mathcal{IF}_{i}^{(2)}\left( y_{i},R_{n,\beta }\left( \boldsymbol{T}%
_{n,\beta }(\underline{F}_{\boldsymbol{\theta }_{2}})\right) \right) \\
& =\mathcal{IF}_{i}^{T}\left( y_{i},\boldsymbol{\lambda }_{n,\beta }\left( 
\boldsymbol{T}_{n,\beta }(\underline{F}_{\boldsymbol{\theta }_{2}})\right)
\right) \left[ \boldsymbol{\Sigma }_{21,\beta }(\boldsymbol{\theta }_{0})%
\boldsymbol{K}_{\beta }(\boldsymbol{\theta }_{0})\boldsymbol{\Sigma }%
_{21,\beta }^{T}(\boldsymbol{\theta }_{0})\right] ^{-1}\mathcal{IF}%
_{i}\left( y_{i},\boldsymbol{\lambda }_{n,\beta }\left( \boldsymbol{T}%
_{n,\beta }(\underline{F}_{\boldsymbol{\theta }_{2}})\right) \right) \\
& =\boldsymbol{d}_{i,\beta }^{T}(\boldsymbol{\theta }_{0})\left[ \boldsymbol{%
M}^{T}(\boldsymbol{\theta }_{0})\boldsymbol{J}_{\beta }^{-1}(\boldsymbol{%
\theta }_{0})\boldsymbol{K}_{\beta }(\boldsymbol{\theta }_{0})\boldsymbol{J}%
_{\beta }^{-1}(\boldsymbol{\theta }_{0})\boldsymbol{M}(\boldsymbol{\theta }%
_{0})\right] ^{-1}\boldsymbol{d}_{i,\beta }(\boldsymbol{\theta }_{0}),
\end{align*}

\section{Computational considerations for obtaining the Breusch-Pagan and Koenker's $\protect\beta $-score test-statistics}
In this section we provide some details about the followed methodology for computation, taking in consideration some mathematical tools as well as some details about the {\texttt{R}} code.
\subsection{Case: $\protect\beta =0$}

\subsubsection{MLE of $\boldsymbol{\protect\theta }=(\boldsymbol{\protect%
		\beta },\protect\sigma ^{2})^{T}$, $\protect\widetilde{\boldsymbol{\protect%
			\theta }}=(\protect\widetilde{\boldsymbol{\protect\beta }},\protect%
	\widetilde{\protect\sigma }^{2})^{T}$}

Let us consider the Cholesky decomposition of $\mathbb{X}_{n}^{T}\mathbb{X}%
_{n}$, i.e. $\mathbb{X}_{n}^{T}\mathbb{X=L}_{n}\mathbb{L}_{n}^{T}$ where $%
\mathbb{L}_{n}$ is an lower triangular matrix . Considering the system of
equations%
\begin{align}
	& \mathbb{X}_{n}^{T}\mathbb{X}_{n}\widetilde{\boldsymbol{\beta }}=\mathbb{X}%
	_{n}^{T}\boldsymbol{Y},  \label{eq.o} \\
	& \mathbb{L}_{n}\mathbb{L}_{n}^{T}\widetilde{\boldsymbol{\beta }}=\mathbb{X}%
	_{n}^{T}\boldsymbol{Y},  \notag
\end{align}%
to be solved in two steps%
\begin{equation*}
	\left\{ 
	\begin{array}{c}
		\mathbb{L}_{n}\widetilde{\widetilde{\boldsymbol{\beta }}}=\mathbb{X}_{n}^{T}%
		\boldsymbol{Y}, \\ 
		\mathbb{L}_{n}^{T}\widetilde{\boldsymbol{\beta }}=\widetilde{\widetilde{%
				\boldsymbol{\beta }}},%
	\end{array}%
	\right.
\end{equation*}%
being the first one forward substitution and the second one backward
substitution.

\subsubsection{Classical Breusch Pagan and Koenker's test-statistics}

Let us consider $v_{i}(\widetilde{\boldsymbol{%
		\theta }})=g_{i}(\widetilde{\boldsymbol{\theta }})-1$, $g_{i}(\widetilde{%
	\boldsymbol{\theta }})=\frac{\varepsilon _{i}^{2}(\widetilde{\boldsymbol{%
			\beta }})}{\widetilde{\sigma }^{2}}$, $i=1,\ldots ,n$ and $\widetilde{\sigma 
}^{2}=\frac{1}{n}\tsum_{i=1}^{n}\varepsilon _{i}^{2}(\widetilde{\boldsymbol{%
		\beta }})$. We must derive an ANOVA of%
\begin{align*}
	\boldsymbol{v}(\widetilde{\boldsymbol{\theta }})& =\mathbb{\breve{Z}}_{n}%
	\boldsymbol{\gamma }+\boldsymbol{\varepsilon } \\
	& =\boldsymbol{1}_{n}\gamma _{0}+\mathbb{Z}_{n}\boldsymbol{\alpha }+%
	\boldsymbol{\varepsilon },
\end{align*}%
where%
\begin{align*}
	\boldsymbol{v}(\widetilde{\boldsymbol{\theta }})& =\boldsymbol{g}(\widetilde{%
		\boldsymbol{\theta }})-\boldsymbol{1}_{n}, \\
	\boldsymbol{g}(\widetilde{\boldsymbol{\theta }})& =\frac{\boldsymbol{%
			\varepsilon }(\widetilde{\boldsymbol{\beta }})\odot \boldsymbol{\varepsilon }%
		(\widetilde{\boldsymbol{\beta }})}{\widetilde{\sigma }^{2}}, \\
	\boldsymbol{\varepsilon }(\widetilde{\boldsymbol{\beta }})& =\boldsymbol{Y}-%
	\mathbb{X}_{n}\widetilde{\boldsymbol{\beta }},
\end{align*}%
$\odot $ is the element wise product of two vectors, $\widetilde{\sigma }^{2}
$ is the solution of $\boldsymbol{v}^{T}(\widetilde{\boldsymbol{\theta }})%
\boldsymbol{1}_{n}=0$. Then%
\begin{align*}
	R_{n}(\widetilde{\boldsymbol{\theta }})& =ESS(\widetilde{\boldsymbol{\theta }%
	})/2, \\
	Q_{n}(\widetilde{\boldsymbol{\theta }})& =\frac{ESS(\widetilde{\boldsymbol{%
				\theta }})}{TSS(\widetilde{\boldsymbol{\theta }})/n}=nR_{OLS}^{2}(\widetilde{%
		\boldsymbol{\theta }}),
\end{align*}%
are respectively the classical Breusch Pagan and Koenker's test-statistics,
which are asymptotically a chi-square random variable with $r$ degrees of
freedom (number of columns of full rank matrix $\mathbb{Z}_{n}$). Under
non-normal errors, the same asymptotic distribution of $Q_{n}(\widetilde{%
	\boldsymbol{\theta }})$ remains being valid but $TSS(\widetilde{\boldsymbol{%
		\theta }})/n$ does not need to tend to $2$ almost surely. The value of $TSS(%
\widetilde{\boldsymbol{\theta }})/n$ is an estimation of the Kurtosis minus $%
1$ of the errors distribution.

\subsection{Case: $\protect\beta >0$}

\subsubsection{Minimum DPD estimation of $\boldsymbol{\protect\theta }=(%
	\boldsymbol{\protect\beta },\protect\sigma ^{2})^{T}$, $\protect\widetilde{%
		\boldsymbol{\protect\theta }}_{\protect\beta }=(\protect\widetilde{%
		\boldsymbol{\protect\beta }}_{\protect\beta },\protect\widetilde{\protect%
		\sigma }_{\protect\beta }^{2})^{T}$}

The solution of%
\begin{equation}
	\mathbb{X}_{n}^{T}\boldsymbol{D}_{\beta }(\widetilde{\boldsymbol{\theta }}%
	_{\beta })\mathbb{X}_{n}\widetilde{\boldsymbol{\beta }}=\mathbb{X}_{n}^{T}%
	\boldsymbol{D}_{\beta }(\widetilde{\boldsymbol{\theta }}_{\beta })%
	\boldsymbol{Y},  \label{eq.w}
\end{equation}%
where%
\begin{equation*}
	\boldsymbol{D}_{\beta }(\widetilde{\boldsymbol{\theta }}_{\beta })=\mathrm{%
		diag}(\exp (-\tfrac{\beta }{2}\boldsymbol{g}(\widetilde{\boldsymbol{\theta }}%
	_{\beta }))),
\end{equation*}%
will be recursively obtained. (\ref{eq.w}) will be solved in the same way as
(\ref{eq.o}), taking%
\begin{align*}
	\mathbb{X}_{n,\beta }(\widetilde{\boldsymbol{\theta }}_{\beta })& =%
	\boldsymbol{D}_{\beta }^{1/2}(\widetilde{\boldsymbol{\theta }}_{\beta })%
	\mathbb{X}_{n}, \\
	\boldsymbol{Y}(\widetilde{\boldsymbol{\theta }}_{\beta })& =\boldsymbol{D}%
	_{\beta }^{1/2}(\widetilde{\boldsymbol{\theta }}_{\beta })\boldsymbol{Y},
\end{align*}%
where%
\begin{equation*}
	\boldsymbol{D}_{\beta }^{1/2}(\widetilde{\boldsymbol{\theta }}_{\beta })=%
	\mathrm{diag}(\exp (-\tfrac{\beta }{4}\boldsymbol{g}(\widetilde{\boldsymbol{%
			\theta }}_{\beta }))),
\end{equation*}%
i.e. (\ref{eq.o}) becomes%
\begin{equation*}
	\mathbb{X}_{n,\beta }^{T}(\widetilde{\boldsymbol{\theta }}_{\beta })\mathbb{X%
	}_{n,\beta }(\widetilde{\boldsymbol{\theta }}_{\beta })\widetilde{%
		\boldsymbol{\beta }}=\mathbb{X}_{n}^{T}(\widetilde{\boldsymbol{\theta }}%
	_{\beta })\boldsymbol{Y}(\widetilde{\boldsymbol{\theta }}_{\beta }).
\end{equation*}
In each iteration $\sigma ^{2}$ is updated solving $\boldsymbol{v}^{T}(%
\widetilde{\boldsymbol{\theta }}_{\beta })\boldsymbol{1}_{n}=0$ in terms of $%
\sigma ^{2}$, i.e. $f(\sigma ^{2})=0$, where%
\begin{equation}
	f(\sigma ^{2})=\frac{1}{n}\tsum_{i=1}^{n}\exp \{-\beta g_{i}(\sigma
	^{2})\}(g_{i}(\sigma ^{2})-1)+\frac{\beta }{(\beta +1)^{\frac{3}{2}}}.
	\label{AppEq1}
\end{equation}%
will be updated.

The $0$-th iteration is just to get $\widetilde{\boldsymbol{\theta }}=(%
\widetilde{\boldsymbol{\beta }},\widetilde{\sigma }^{2})^{T}$ or $\widetilde{%
	\boldsymbol{\theta }}_{\beta =0}=(\widetilde{\boldsymbol{\beta }}_{\beta =0},%
\widetilde{\sigma }_{\beta =0}^{2})^{T}$. The $i$-th iteration consists on
updating $\mathbb{X}_{n,\beta }(\widetilde{\boldsymbol{\theta }}_{\beta })$
and $\boldsymbol{Y}(\widetilde{\boldsymbol{\theta }}_{\beta })$, solving the
system of equations (based on $\widetilde{\boldsymbol{\beta }}_{\beta
	,(i-1)} $, $\widetilde{\sigma }_{\beta ,(i-1)}^{2}$) to get $\widetilde{%
	\boldsymbol{\beta }}_{\beta ,(i)}$, update of $\widetilde{\sigma }_{\beta
	,(i)}^{2}$ solving a non-linear equation with a unique variable and checking
whether the norm of $\mathbb{X}_{n}^{T}(\widetilde{\boldsymbol{\theta }}%
_{\beta })\left( \boldsymbol{Y}(\widetilde{\boldsymbol{\theta }}_{\beta })%
\mathbb{-X}_{n,\beta }(\widetilde{\boldsymbol{\theta }}_{\beta })\widetilde{%
	\boldsymbol{\beta }}\right) $ is less than a pre-specified tolerance level.
If so, the sequence of iterations stops, if not the sequence goes ahead. The
system of linear equations can be efficiently solved using the base R
function \texttt{solve()} to solve linear equations, once the Cholesky
decomposition is properly suited according to the indications given
previously:

\texttt{library(MASS)}

\texttt{XtX \TEXTsymbol{<}- crossprod(X.tilde)}

\texttt{if (is.singular.matrix(XtX,tol=1e-08))\{}

\texttt{print("singular cross-product matrix")}

\texttt{break}

\texttt{\} }

\texttt{library(matrixcalc)}

\texttt{L \TEXTsymbol{<}- t(chol(XtX))}

\texttt{beta.hat \TEXTsymbol{<}- solve(t(L), solve(L, t(X.tilde) \%*\%
	Y.tilde))}

The non-linear equation can be efficiently solved utilizing the \texttt{R}
package \texttt{nleqslv} (Klein and Sporleder (2020), version 0.8-2,
installed on 2020-12-01.), as it enables the attainment of a solution within
a minimal number of iterations. It needs providing (\ref{AppEq1}) and 
\begin{equation*}
	\frac{\partial }{\partial \sigma ^{2}}f(\sigma ^{2})=\frac{\beta }{n\sigma
		^{2}}\sum_{i=1}^{n}\exp \{-\beta g_{i}(\sigma ^{2})\}g_{i}^{2}(\sigma ^{2})-%
	\frac{\beta +1}{n\sigma ^{2}}\sum_{i=1}^{n}\exp \{-\beta g_{i}(\sigma
	^{2})\}g_{i}(\sigma ^{2}),
\end{equation*}%
where%
\begin{align*}
	g_{i}(\sigma ^{2})& =\frac{r_{i}^{2}(\boldsymbol{\beta })}{\sigma ^{2}}, \\
	\frac{\partial }{\partial \sigma ^{2}}g_{i}(\sigma ^{2})& =\left( -\frac{1}{%
		\sigma ^{2}}\right) \frac{r_{i}^{2}(\boldsymbol{\beta })}{\sigma ^{2}}%
	=\left( -\frac{1}{\sigma ^{2}}\right) g_{i}(\sigma ^{2}),
\end{align*}%
with the following code:
\begin{verbatim}
	get.eq <- function(s2) {
		g <- (res*res)/s2
		eg2 <- exp(-beta*g/2.0)
		v <- (eg2) * (g-1.0) #(D.half.beta %*% D.half.beta) %*% (g-1.0)
		v <- mean(v)+beta/((beta+1.0)^1.5)
		return(v)
	}
	get.grad.eq <- function(s2) {
		g <- (res*res)/s2
		eg2 <- exp(-beta*g/2.0)
		dv <- beta * (eg2) * (g*g) #(D.half.beta %*% D.half.beta) %*% (g-1.0)
		dv <- dv - (beta+1.0) * (eg2) * g
		dv <- mean(dv)/s2
		return(dv)
	}  
	library(nleqslv)
	initialG <- c(sigma2.hat)
	sp <- list(trace=1,allowSingular=TRUE,btol = 1.e-5)
	zero <- nleqslv(initialG, get.eq, jac=get.grad.eq, method="Newton",control=sp)
	sigma2.hat <- zero$x
\end{verbatim}

\subsubsection{Breusch-Pagan and Koenker's $\protect\beta $-score
	test-statistic}

In a similar way done for the classical ones, we get%
\begin{align*}
	R_{n}(\widetilde{\boldsymbol{\theta }}_{\beta })& =\frac{ESS(\widetilde{%
			\boldsymbol{\theta }}_{\beta })}{\tfrac{2(2\beta ^{2}+1)}{(2\beta +1)^{5/2}}-%
		\frac{\beta ^{2}}{(\beta +1)^{3}}}, \\
	Q_{n}(\widetilde{\boldsymbol{\theta }}_{\beta })& =\frac{ESS(\widetilde{%
			\boldsymbol{\theta }}_{\beta })}{TSS(\widetilde{\boldsymbol{\theta }}_{\beta
		})/n}=nR_{OLS}^{2}(\widetilde{\boldsymbol{\theta }}_{\beta }).
\end{align*}%
Therefore, once we have calculated $\widetilde{\boldsymbol{\theta }}_{\beta }
$, the scheme to be followed is exactly the same as the one we have followed for
the classical tests and it is easily calculated from the well-known ANOVA
for OLS, hence the {\texttt{ln()}} function and related ones of {\texttt{R}} can be used.
\section*{Appendix References}
\vspace{-30pt} 
\renewcommand{\refname}{}


\begin{thebibliography}{99}
\bibitem{Alih} Alih, E. and Ong, H.C. (2015). An outlier-resistant test for
heteroskedasticity in linear models. \emph{Journal of Applied Statistics}, 
\textbf{42}, 1617--1634.

\bibitem{Aitchin} Aitchison, J. and Silvey, D.S. (1958). Maximum-Likelihood
Estimation of Parameters Subject to Restraints. \emph{The Annals of
Mathematical Statistics}, \textbf{29}, 813--828.

\bibitem{Athe} Athreya, K.B. and Lahiri, S.N. (2006). \emph{Measure Theory
and Probability Theory}. Springer-Verlag.

\bibitem{Basu0} Basu, A., Ghosh, A., Martin, N., and Pardo, L. (2022). A
Robust Generalization of the Rao Test. \emph{Journal of Business \& Economic
Statistics}, \textbf{40}, 868--879.

\bibitem{Basu1} Basu, A, Harris, I.R., Hjort, N.L. and Jones, M.C. (1998).
Robust and efficient estimation by minimising a density power divergence. 
\emph{Biometrika}, \textbf{85}, 549--559.

\bibitem{Basu2} Basu, A., Mandal, A., Martin, N., and Pardo, L. (2017).
Testing Composite Hypothesis Based on the Density Power Divergence. \emph{%
Sankhya B}. \textbf{80}, 222--262.

\bibitem{Basu3} Basu, A., Shioya, H., and Park C. (2011). \emph{Statistical
Inference: The Minimum Distance Approach}. Boca Raton, CRC Press

\bibitem{Beran} Beran, R. J. (1977). M\'{\i}nimum Hellinger distance
estimates for parametric models. \emph{The annals of Statistics}, \textbf{5}%
, 445--463.

\bibitem{Berenguer} Berenguer-Rico, V. and Wilms, I. (2021).
Heteroscedasticity testing after outlier removal. \emph{Econometric Reviews}%
, \textbf{40}, 51--85.

\bibitem{Bickel} Bickel, P.J. (1978). Using residuals robustly I: Tests for
heteroscedasticity and non-linearity. \emph{Annals of Statistics}, 6,
266--291.

\bibitem{Br0} Breusch, T. S. (1978). Testing for Autocorrelation in Dynamic
Linear Models. \emph{Australian EconomicPapers}, \textbf{17}, 334--355.

\bibitem{Breusch} Breusch, T.S. and Pagan, A.R. (1979). A simple test for
heteroscedasticity and random coefficient variation. \emph{Econometrica}, 
\textbf{47}, 1287--1294.

\bibitem{Breusch2} Breusch, T.S., and Pagan, A.R. (1980). The Lagrange
multiplier test and its applications to model specification in econometrics. 
\emph{The Review of Economic Studies}, 47, 239--253.

\bibitem{Boos} Boos, D. D. (1992). On Generalized Score Tests. \emph{The
American Statistician}, \textbf{46}, 327--333.

\bibitem{Carrol} Carrol, R.J. and Ruppert, D. (1988). \emph{Transformation
and Weighting in Regression}. CRC Press.

\bibitem{Cook} Cook, R. D. and Weisberg, S. (1983). Diagnostics for
heteroscedasticity in regression. \emph{Biometrika}, \textbf{70}, 1--10.

\bibitem{Cressie} Cressie, N. and Read, T.R.C. (1984). Multinomial
goodness-of-fit tests. \emph{Journal of the Royal Statistical Society.
Series B}, \textbf{46}, 440--464.

\bibitem{Sciszar} Csisz\'{a}r, I. (1967). Information-type measures of
difference of probability distributions and indirect observations. \emph{%
Studia Scientiarum Mathematicarum Hungarica}, \textbf{2}, 299--318.

\bibitem{Ghosh1} Ghosh, A. and Basu, A. (2013). Robust estimation for
independent non-homogeneous observations using density power divergence with
applications to linear regression. \emph{Electronic Journal of Statistics}, 
\textbf{7}, 2420--2456.

\bibitem{Ghosh2} Ghosh, A., Mandal, A., Mart\'{\i}n, N. and Pardo, L.
(2016). Influence analysis of robust Wald-type tests. \emph{Journal of
Multivariate Analysis}, \textbf{147}, 102--126.

\bibitem{Godfrey00} Godfrey, L.G. (1978a). Testing for multiplicative heteroscedasticity.
\emph{Journal of Econometrics}, \textbf{8}, 227-236.

\bibitem{Godfrey0} Godfrey, L.G. (1978b). Testing Against General
Autoregressive and Moving Average Error Models when the Regressors Include
Lagged Dependent Variables. \textit{Econometrica}, \textbf{46}, 1293--1301.

\bibitem{Godfrey} Godfrey, L. G. (1989). \emph{Misspecification Tests in
Econometrics: The Lagrange Multiplier Principle and Other Approaches}.
Cambridge University Press, New York.

\bibitem{Halunga} Halunga, A.G., Orme, C.D, Yamagata, T. (2017). A
heteroskedasticity robust Breusch-Pagan test for Contemporaneous correlation
in dynamic panel data models. \emph{Journal of Econometrics}. 198, 209--230.

\bibitem{Hampel0} Hampel, F.R. (1968). \emph{Contribution to the theory of
robust estimation}. Ph.D. Thesis, University of California, Berkeley.

\bibitem{Hampel} Hampel, F.R. (1974). The influence curve and its role in
robust estimation. \emph{Journal of the American Statistical Association}, 
\textbf{69}, 383--393.

\bibitem{Hampel2} Hampel, F.R., E.M. Ronchetti, P.J. Rousseeuw, and W.A.
Stahel (1986). \emph{Robust Statistics: The Approach Based on lnfluence
Functions}. Wiley, New York.

\bibitem{Hannan} Hannan, E.J. (1956). The Asymptotic Powers of Certain Tests
Based on Multiple Correlations. \emph{Journal of the Royal Statistical
Society. Series B}, \textbf{18}, 227--233.

\bibitem{Honda} Honda, Y. (1988). A size correction to the Lagrange
multiplier test for heteroskedasticity. \emph{Journal of Econometrics}, 
\textbf{38}, 375--386.

\bibitem{Kalina} Kalina, J. (2011) Testing heteroscedasticity in robust
regression. \emph{Research Journal of Economics Business and ICT}, \textbf{4}%
, 25--28.

\bibitem{Kim1} Kim, B. and Lee, S. (2013). Robust estimation for copula Parameter in SCOMDY
models. Journal of Time Series Analysis, \textbf{34}, 302--314.

\bibitem{Kim2} Kim, B. and Lee, S. (2017). Robust estimation for zero-inflated Poisson
autoregressive models based on density power divergence. Journal of
Statistical Computation and Simulation, \textbf{87}, 2981--2996.

\bibitem{Kim2} Kim, B. and Lee, S. (2020). Robust estimation for general integer-valued time
series models. Annals of the Institute of Statistical Mathematics, 72(6), 1371-1396.

\bibitem{Kim2} Kim, B. (2018). Robust maximum entropy test for GARCH models based on a
minimum density power divergence estimator. Economics Letters, 162, 93--97.

\bibitem{Lee} Lee, S. and Song, J. (2009). Minimum density power divergence estimator for
GARCH models. TEST 18, 316--341.

\bibitem{Koenker1} Koenker, R. (1981). A note on Studentizing a test for
heteroscedasticity. \emph{Journal of Econometrics}, \textbf{17}, 107--112.

\bibitem{Koenker2} Koenker, R. and Bassett, G. (1982). Robust tests for
heteroscedasticity based on regression quantiles. \emph{Econometrica}, 
\textbf{50}, 43--61.

\bibitem{Koenker3} Koenker, R., and Bassett, G. (1982). Tests of Linear
Hypotheses and $l_{1}$ Estimation. \emph{Econometrica}, \textbf{50},
1577--1583.

\bibitem{krasker} Krasker, W.S. and Welsch, R.E. (1982). Efficient bounded
influence regression estimation. \emph{Journal of the American Statistical
Association}, \textbf{77}, 595--604.

\bibitem{Kolmogorov} Kolmogorov, A. N. (1933). Sulla determinazione empirica
di una legge di distribuzione. \emph{Giornale dell'Istituto Italiano degli
Attuari}, \textbf{4}, 83-91.

\bibitem{Kmenta} Kmenta, J. (1986). \emph{Elements of Econometrics}.
Macmillan.

\bibitem{Lyon} Lyon, J. D. and Tsai, C.L. (1996). A Comparison of Tests for
Heteroscedasticity. \emph{Journal of the Royal Statistical Society. Series D
(The Statistician)}, \textbf{45}, 337--349.

\bibitem{Lu} Lu, T. T. and Shiou, S. H. (2002). Inverses of $2\times 2$
block matrices. \emph{Computers \& Mathematics with Applications}, \textbf{43%
}, 119--129.

\bibitem{Martin} Mart\'{\i}n, N. (2021). Rao's Score Tests on Correlation
Matrices. arXiv, Statistics Theory (math.ST),
\texttt{\href{https://doi.org/10.48550/arXiv.2012.14238}%
	{https://doi.org/10.48550/arXiv.2012.14238}}.

\bibitem{Mahalanobis} Mahalanobis, P. C. (1936). On the generalized distance
in statistics. \emph{Proceedings of the National Institute of Science of
India}, \textbf{2}, 49--55.

\bibitem{Markatou} Markatou, M., Karlis, D., and Ding, Y. (2021).
Distance-Based Statistical Inference. \emph{Annual Review of Statistics and
Its Application}, \textbf{8}, 301--327.

\bibitem{Nuttal} A. Nuttall (1975). Some integrals involving the $Q_{M}$\
function. \emph{IEEE Transactions on Information Theory}, \textbf{21},
95--96.

\bibitem{Pardo} Pardo, L. (2006). Statistical inference based on divergence
measures.\emph{\ Chapman and Hall/CRC}.

\bibitem{Rao} Rao, C. R., Toutenburg, H., Shalabh and Heumann, C. (2008). 
\emph{Linear Models and Generalizations. Least Squares and Alternatives}.
Springer.

\bibitem{Salivian} Salibian-Barrera, M., Van Aelst, S. and Yohai, V.J.
(2016). Robust tests for linear regression models based on $\tau $%
-estimates. \emph{Computational Statistics \& Data Analysis}, \textbf{93},
436--455.

\bibitem{Silvey} Silvey, S.D. (1959). The Lagrangian Multiplier Test. \emph{%
The Annals of Mathematical Statistics}, \textbf{30}, 389--407.

\bibitem{Waldman} Waldman, D.M. (1983). A note on algebraic equivalence of White's test and a
variation of the Godfrey/Breusch-Pagan test for heteroscedasticity.
\emph{Economics Letters}, 13, 197--200.

\bibitem{White} White. H. (1980). A heteroscedasticity-consistent covariance matrix estimator
and a direct test for heteroscedasticity. \emph{Econometrica}, 48, 817-838. 

\end{thebibliography}

\begin{thebibliography}{9}
	\bibitem{Davidson} Davidson, R., and MacKinnon, J. G. (2021). \emph{%
		Estimation and Inference in Econometrics}. Oxford University Press, New York.
	\bibitem{klein} Klein, S. and Sporleder, C. (2020). nleqslv: Solve Nonlinear
	Equations and systems of Equations. R package version 0.8-2. \texttt{\href{https://CRAN.R-project.org/package=nleqslv}%
		{https://CRAN.R-project.org/package=nleqslv}}
\end{thebibliography}
\end{document}